\documentclass[12pt]{amsart}
\usepackage[T1]{fontenc}
\usepackage[english]{babel}
\usepackage{mathtools}
\usepackage{graphicx}
\usepackage{amssymb,bbm}
\usepackage{tikz}
\usetikzlibrary{arrows}
\usepackage{hyperref}
\usepackage{dsfont,soul}
\usepackage{enumitem}
\usepackage{multirow}
\usepackage{subfigure}
\graphicspath{{graphics/}}
\usepackage{geometry}
\setlist[enumerate]{labelsep=*, leftmargin=2pc,
topsep=1ex plus0.5ex minus0.2ex,
itemsep=1ex plus0.5ex minus0.2ex,
font=\rmfamily,
font=\upshape}
\setlist[itemize]{labelsep=*, leftmargin=1pc,
topsep=1ex plus0.5ex minus0.2ex,
itemsep=1ex plus0.5ex minus0.2ex,
font=\rmfamily,
font=\upshape}
\newtheorem{thm}{Theorem}[section]

\newtheorem{lem}[thm]{Lemma}

\theoremstyle{definition}

\newtheorem{exa}[thm]{Example}
\newtheorem{rem}[thm]{Remark}
\numberwithin{equation}{section}
\newcommand{\bC}{\mathbb C}
\newcommand{\bE}{\mathbb E}
\newcommand{\bF}{\mathbb F}
\newcommand{\bN}{\mathbb N}
\newcommand{\bP}{\mathbb P}
\newcommand{\bR}{\mathbb R}
\newcommand{\cA}{\mathcal A}
\newcommand{\cE}{\mathcal E}

\newcommand{\cL}{\mathcal L}
\newcommand{\cM}{\mathcal M}
\newcommand{\cN}{\mathcal N}

\newcommand{\her}{M^{\rm h}}

\newcommand{\id}{\mathds{1}} 
\DeclareMathOperator{\ii}{i} 
\DeclareMathOperator*{\argmax}{argmax}
\DeclareMathOperator{\conv}{conv}
\DeclareMathOperator{\diag}{diag}
\DeclareMathOperator{\tr}{tr}
\begin{document}
\selectlanguage{english}
\title{Classification of joint numerical ranges of three
hermitian matrices of size three}
\author{Konrad Szyma\'nski, Stephan Weis, and Karol {\.Z}yczkowski}
\begin{abstract}
The joint numerical range $W(F)$ of three hermitian $3$-by-$3$ matrices 
$F=(F_1,F_2,F_3)$ is a convex and compact subset in $\bR^3$. We show 
that $W(F)$ is generically a three-dimensional oval.
Assuming $\dim(W(F))=3$, every one- or two-dimensional face of $W(F)$ is 
a segment or a filled ellipse. We prove that only ten configurations of
these segments and ellipses are possible. We identify a triple $F$ for
each class and illustrate $W(F)$ using random matrices and dual varieties.
\end{abstract}
%
%
%
\date{February 11, 2018}
\maketitle
\thispagestyle{empty}
\pagestyle{myheadings}
\markleft{\hfill Classification of joint numerical ranges \hfill}
\markright{\hfill K.\ Szyma\'nski et al.\ \hfill}
%
%
\section{Introduction}
\par
We denote the space of complex $d$-by-$d$ matrices by $M_d$, the real subspace 
of hermitian matrices by $\her_d:=\{a\in M_d\mid a^*=a\}$, and the identity 
matrix by $\id_d\in M_d$. We write 
$\langle x,y\rangle:=\overline{x_1}y_1+\cdots+\overline{x_d}y_d$, 
$x,y\in\bC^d$ for the inner product on $\bC^d$. 
For $d,n\in\bN=\{1,2,3,\ldots\}$, let $F:=(F_1,\ldots,F_n)\in(\her_d)^n$ be
an $n$-tuple of hermitian $d$-by-$d$ matrices. The {\it joint numerical range} 
of $F$ is 
\[
W(F):=\{(\langle x,F_1x\rangle,\ldots,\langle x,F_nx\rangle)
\mid x\in\bC^d, \langle x,x\rangle=1\}
\subset\bR^n.
\]
For $n=2$, identifying $\bR^2\cong\bC$, the set $W(F_1,F_2)$ is the 
{\em numerical range} 
$\{\langle x,Ax\rangle\mid x\in\bC^d, \langle x,x\rangle=1\}$
of $A:=F_1+\ii F_2$. The numerical range is convex for all $d\in\bN$
by the Toeplitz-Hausdorff theorem \cite{Toeplitz1918,Hausdorff1919}. 
Similarly, for $n=3$ and all $d\geq 3$ the joint numerical range 
$W(F_1,F_2,F_3)$ is convex \cite{Au-YeungPoon1979}. However, $W(F)$ 
is in general not convex for $n\geq 4$, see 
\cite{Polyak1998,LiPoon2000,Gutkin-etal2004}.
\par
Let $d,n\in\bN$ be arbitrary and call $F\in(\her_d)^n$ {\em unitarily reducible} 
if there is a unitary $U\in M_d$ such that the matrices 
$U^*F_1U,\ldots,U^*F_nU$ have a common block diagonal form with two proper blocks. 
Otherwise $F$ is {\em unitarily irreducible}. 
\par
The shape of the numerical range ($n=2$) is well understood. The elliptical 
range theorem \cite{Li1996} states that the numerical range of a $2$-by-$2$
matrix is a singleton, segment, or filled ellipse. Kippenhahn 
\cite{Kippenhahn1951} proved for all $d\in\bN$ and $F\in(\her_d)^2$ that 
$W(F)$ is the convex hull of the boundary generating curve defined in 
Remark~\ref{rem:intro}. He showed for $3$-by-$3$ matrices ($d=3$) that if 
$F$ is unitarily reducible, then $W(F)$ is a singleton, segment, triangle, 
ellipse, or the convex hull of an ellipse and a point outside the ellipse. 
If $F$ is unitarily irreducible, then $W(F)$ is an ellipse, the convex hull 
of a quartic curve (with a flat portion on the boundary), or the convex hull 
of a sextic curve (an oval). Kippenhahn's result for $3$-by-$3$ matrices 
was expressed in terms of matrix invariants and matrix entries of $F_1+\ii F_2$, 
see \cite{Keeler-etal1997,RodmanSpitkovsky2005,Rault-etal2013,SpitkovskyWeis2015}.
The boundary generating curve was also used \cite{ChienNakazato2012} to find
a classification of the numerical range of a $4$-by-$4$ matrix. Another
result \cite{Henrion2010,HeltonSpitkovsky2012} is that a subset $W$ of
$\bC$ is the numerical range of some $d$-by-$d$ matrix if and only if it 
is a translation of the polar of a rigidly convex set of degree less than 
or equal to $d$, see Corollary 3 of \cite{HeltonSpitkovsky2012}. We omit 
the details of this last description as we will not use it.
%
%
%
\begin{figure}
    \centering
    \subfigure[Ex.~\ref{exa:01}, $s=0, e=1$]{
        \includegraphics[width=3.4cm]{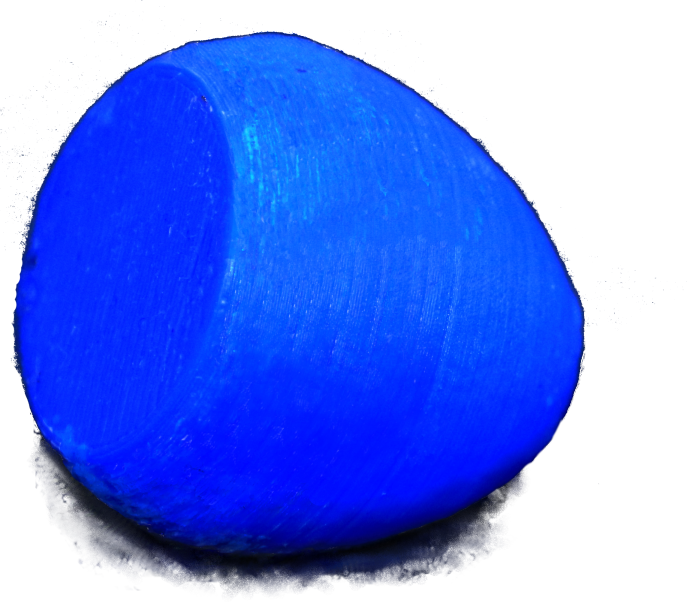}
                \label{fig:print01}
    }
    \subfigure[Ex.~\ref{exa:03}, $s=0, e=3$]{
        \includegraphics[width=3.4cm]{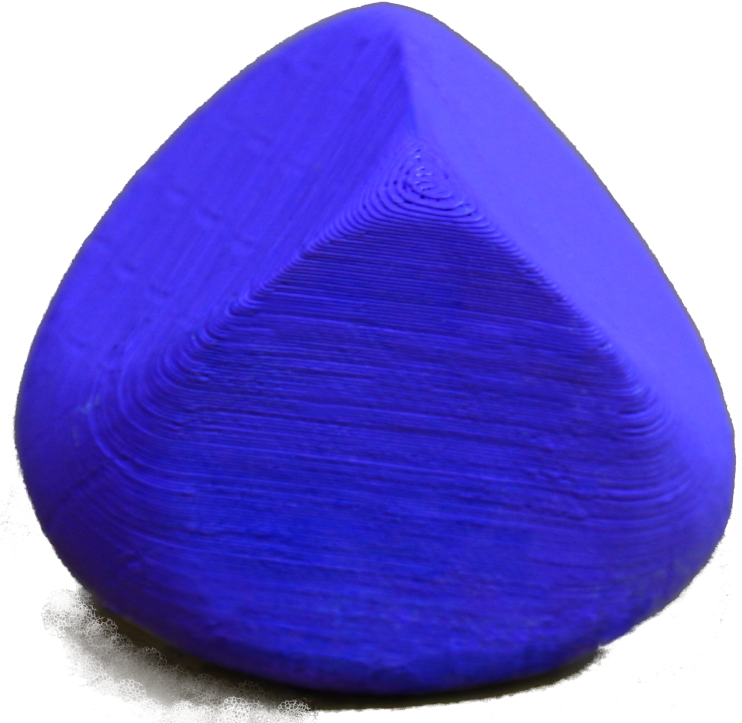}
                 \label{fig:print03}
    }
    \subfigure[Ex.~\ref{exa:04}, $s=0, e=4$]{
        \includegraphics[width=3.4cm]{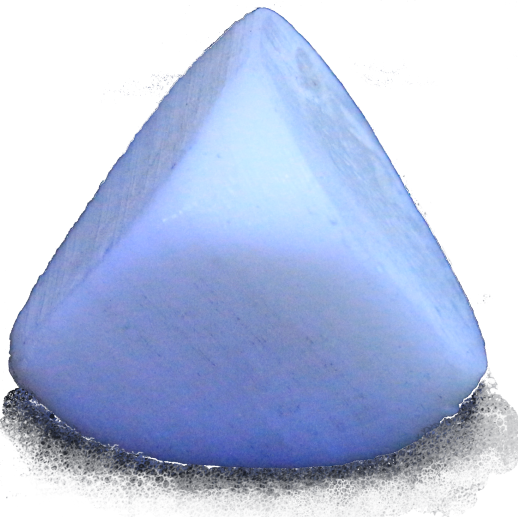}
                 \label{fig:print04}
    }
    \subfigure[Ex.~\ref{exa:10}, $s=1, e=0$]{
        \includegraphics[width=3.4cm]{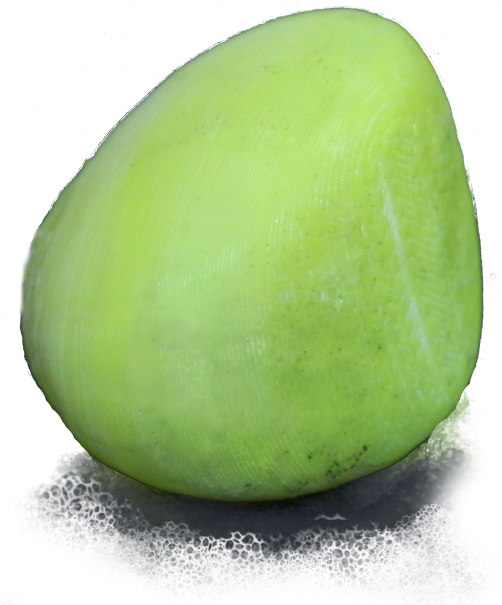}
                 \label{fig:print10}
    }
        \subfigure[Ex.~\ref{exa:11}, $s=1, e=1$]{
        \includegraphics[width=3.4cm]{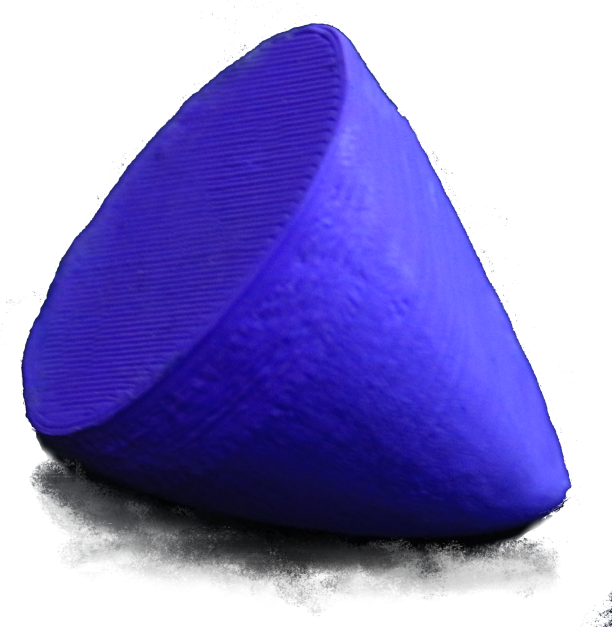}
                 \label{fig:print11}
    }
    \subfigure[Ex.~\ref{exa:12}, $s=1, e=2$]{
        \includegraphics[width=3.4cm]{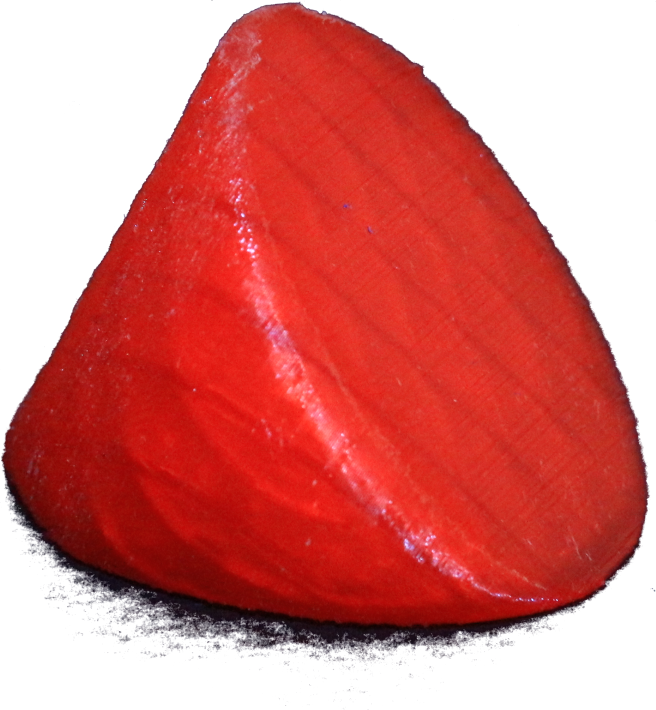}
                 \label{fig:print12}
    }
    \subfigure[$s=e=0$]{
        \includegraphics[width=3.4cm]{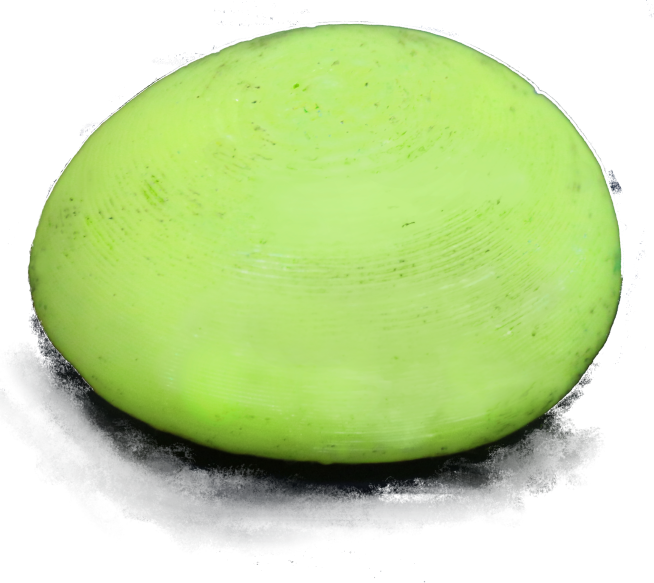}
                 \label{fig:print_generic}
    }
    \subfigure[$s=\infty, e=0$]{
        \includegraphics[width=3.4cm]{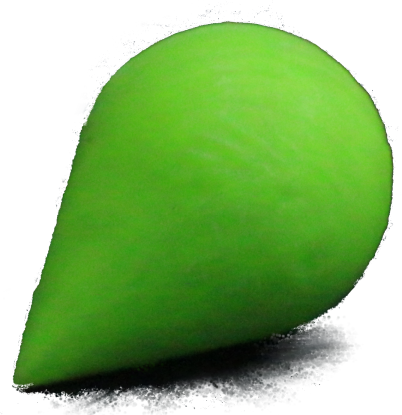}
                 \label{fig:print_leprechaun}
    }
    \subfigure[$s=\infty, e=1$]{
        \includegraphics[width=3.4cm]{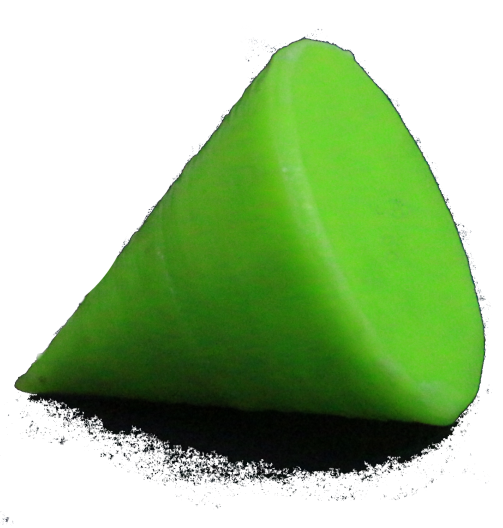}
                \label{fig:print_cone}
    }
    \caption{%
    3D printouts of exemplary joint numerical ranges of triples of hermitian 
    3-by-3 matrices from random density matrices: $s$ denotes the number of 
    segments, $e$ the number of ellipses in the boundary.}
    \label{fig:print}
\end{figure}%
%
%
%
\par
Despite the long history of the problem
\cite{BindingLi1991,KrupnikSpitkovsky2006,ChienNakazato2009,
ChienNakazato2010,Cheung-etal2011}, a classification of the joint numerical
range of triples of matrices ($n=3$) is unknown even in the case $d=n=3$. 
Our motivation to tackle this problem is quantum mechanics, as we explain in 
Section~\ref{sec:qm}. The link to physics is that for arbitrary $d,n\in\bN$ 
and $F\in(\her_d)^n$ the convex hull 
\[
L(F):=\conv(W(F))
\]
of $W(F)$ is a projection (image under a linear map) of the state space
$\cM_d$ of the algebra $M_d$ \cite{DGHMPZ11}. The state space consists of 
$d$-by-$d$ {\em density matrices}, that is positive semi-definite matrices 
of trace one, which represent the states of a quantum system. 
\par
Until further notice let $d=n=3$, where $L(F)=W(F)$ holds. One of us used
random matrices to compute exemplary joint numerical ranges
\cite{Zyczkowski-etal2011}. Photos of their printouts on a 3D-printer are 
depicted in Figure~\ref{fig:print}. The printout shown in 
Figure~\ref{fig:print12} was the starting point of this research. As a 
result we present a simple classification of $W(F)$ in terms of exposed 
faces.  An {\em exposed face} of $W(F)$ is a subset of $W(F)$ which is 
either empty or consists of the maximizers of a linear functional on $W(F)$. 
Lemma~\ref{lem:large-faces-exposed}  shows that the non-empty exposed faces 
of $W(F)$ which are neither singletons nor equal to $W(F)$ are segments or 
filled ellipses. We call them {\em large faces} of $W(F)$ and collect them 
in the set
\begin{align}\label{eq:large-face}
\cL(F)  \,:=\, \{  & \mbox{$G$ is an exposed face of $W(F)$} \mid \\
 & \mbox{$G\neq W(F)$ and $G$ is a segment or a filled ellipse} \}.\nonumber
\end{align} 
Let $e$ (resp.~$s$) denote the number of filled ellipses (resp.\ segments) 
in $\cL(F)$. We recall that a {\em corner point} of $W(F)$ is a point which 
lies on three supporting hyperplanes with linearly independent normal vectors.
\begin{figure}
\includegraphics[height=4cm]{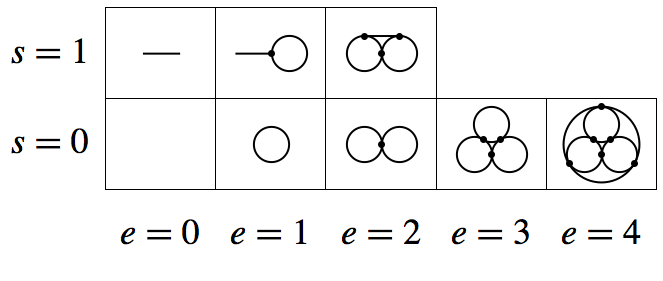}
\caption{\label{fig:classes}%
Possible configurations of large faces of a joint numerical range without 
corner points for $d=n=3$. Circles (resp.~segments) denote large faces 
which are filled ellipses (resp.~segments). Dots denote intersection points 
between large faces.}
\end{figure}
\begin{thm}\label{thm:intro}
Let $F\in(\her_3)^3$. If $W(F)$ has no corner point, then the set 
$\cL(F)$ of large faces of $W(F)$ has one of the eight configurations 
of Figure~\ref{fig:classes}.
\end{thm}
\par\noindent
{\em Proof:}
It is easy to see that large faces intersect mutually (Lemma~\ref{lem:class1}). 
Since $W(F)$ has no corner point, no point lies on three mutually distinct 
large faces (Lemma~\ref{lem:class2}). Hence the union of large faces contains 
an embedded complete graph with one vertex at the centroid of each large face 
(Lemma~\ref{lem:graph-embedding}). Now, a well-known theorem of graph 
embedding \cite{RingelYoungs1968} shows $e+s\leq 4$. We observe that $s=0,1$ 
holds, because for $s\geq 2$ the set $W(F)$ has a corner point 
(Lemma~\ref{lem:exclude-two-segments}). We exclude the case $(e,s)=(3,1)$ by 
noting that for $s\geq 1$ the embedded complete graph has a vertex on a 
segment. Then the graph has vertex degree at most two which implies $e+s\leq 3$.
\hspace*{\fill}$\square$\\
\par
Section~\ref{sec:examples} shows three-dimensional examples of $W(F)$ without 
corner points for all configurations of Figure~\ref{fig:classes}. We are
unaware of earlier examples of
\[
(e,s)\quad = \quad
(1,0),\quad
(2,0),\quad
(3,0),\quad
\mbox{and}\quad
(0,1).
\]
Ovals, where $(e,s)=(0,0)$, are studied in 
\cite{KrupnikSpitkovsky2006}. An example of $(e,s)=(4,0)$ is in 
\cite{Henrion2011}, one of $(e,s)=(1,1)$ is in \cite{ChienNakazato2010}, and 
one of $(e,s)=(2,1)$ is in \cite{Chen-etal2015}.
\par
If $\dim(W(F))=3$ and $W(F)$ has corner points, then 
Lemma~\ref{lem:corner-of-3-by-3} shows that $W(F)$ is the convex hull of 
an ellipsoid and a point outside the ellipsoid, where $(e,s)=(0,\infty)$, 
or the convex hull of an ellipse and a point outside the affine hull of 
the ellipse, where $(e,s)=(1,\infty)$. Examples are depicted in 
Figure~\ref{fig:print_leprechaun} and~\ref{fig:print_cone}.
\par
If $\dim(W(F))=2$ then $e=0$. By projecting to a plane, $W(F)$ 
corresponds to the numerical range of a 3-by-3 matrix. Notice that $W(F)$ 
belongs to one of four classes of 2D objects characterized by the number 
of segments $s=0,1,2,3$. The classification of $W(F)$ in terms of this 
number $s$ is courser than that explained above \cite{Kippenhahn1951}. An 
object with $s=0$ can be an ellipse or the convex hull of a sextic curve.
\begin{rem}[Limits of extreme points] 
Three-dimensional joint numerical ranges of $F\in(\her_3)^3$ solve a 
problem posed in \cite{Weis2012a}. A limit of extreme points 
of $\cM_d$, $d\in\bN$, is again an extreme point and the question is 
whether the analogue holds for projections of $\cM_d$. This doubt is 
dispelled by observing for $(e,s)=(0,1),(1,1),(2,1)$ that any point 
in the relative interior of the segment in $\cL(F)$ is a limit of 
extreme points of $W(F)$ but no extreme point itself, see 
Figure~\ref{fig:2}. The problem was already solved in Example 6 of 
\cite{Chen-etal2015} and discussed in Example~4.2 of \cite{Rodman-etal2016} 
with an example of $(e,s)=(2,1)$. A simpler example, with larger 
matrices, is 
\[
F_1=\left(\begin{smallmatrix}
1 &  0 & 0 & 0 \\
0 & -1 & 0 & 0 \\
0 &  0 & 1 & 0 \\
0 &  0 & 0 & 1 
\end{smallmatrix}\right),
\qquad
F_2=\left(\begin{smallmatrix}
0 & 1 & 0 & 0 \\
1 & 0 & 0 & 0 \\
0 & 0 & 0 & 0 \\
0 & 0 & 0 & 0 
\end{smallmatrix}\right),
\qquad
F_3=\left(\begin{smallmatrix}
0 & 0 &  0 & 0 \\
0 & 0 &  0 & 0 \\
0 & 0 & -1 & 0 \\
0 & 0 &  0 & 1 
\end{smallmatrix}\right),
\]
where $W(F_1,F_2,F_3)$ is the convex hull of the union of the unit disk 
in the $x$-$y$-plane with the points $(1,0,\pm1)$.
\end{rem}
\begin{rem}[Real varieties] 
\label{rem:intro}
For arbitrary $d,n\in\bN$ we consider the hypersurface in the complex 
projective space $\bP^n$, defined as the zero locus
\[\textstyle
S_F:=
\{(u_0:\cdots:u_n)\in\bP^n\mid \det(u_0\id+u_1F_1+\cdots+u_nF_n)=0\}.
\]
An analysis of singularities of $S_F$ for $d=n=3$ shows that $W(F)$ has 
at most four large faces which are ellipses \cite{ChienNakazato2009}. 
This estimate also follows from our classification. The {\em dual variety} 
$S_F^*\subset\bP^n{}^*$ is the complex projective variety which is the 
closure of the set of tangent hyperplanes of $S_F$ at smooth points 
\cite{Fischer2001,Harris1995,Gelfand-etal1994}. The 
{\em boundary generating hypersurface} \cite{ChienNakazato2010}
of $F$ is the real affine part of the dual variety,
\[
S_F^*(\bR)
:=\{(x_1,\ldots,x_n)\in\bR^n\mid(1:x_1:\cdots:x_n)\in S_F^*\}.
\]
For $n=2$, the variety $S_F^*(\bR)$ is called 
{\em boundary generating curve}, and Kippenhahn \cite{Kippenhahn1951}
showed that the convex hull of $S_F^*(\bR)$ is the numerical range
$W(F_1,F_2)$. A more detailed proof is given in \cite{ChienNakazato2010}. 
For $d=n=3$, Chien and Nakazato \cite{ChienNakazato2010} discovered 
that $S_F^*(\bR)$ can contain (unbounded) lines, so the analogue of 
the Kippenhahn assertion is wrong for $n\geq 3$. We will see examples 
of such lines in Section~\ref{sec:examples}.
\end{rem}
\par
Section~\ref{sec:boundary} studies exposed faces. One result is that the 
joint numerical range of $F\in(\her_3)^3$ is generically an {\em oval}, 
that is a compact strictly convex set with interior points and smooth 
boundary. More generally, Theorem~\ref{thm:generic} shows that 
$L(F)=\conv(W(F))$ is generically an oval for all $d\geq 2$ and 
$n\leq 3$, using the von Neumann-Wigner {\em non-crossing rule} 
\cite{NeumannWigner1929,Gutkin-etal2004} and results about normal cones 
developed in Section~\ref{sec:lattices}. Using the {\em crossing rule} 
\cite{Friedland-etal1984}, Lemma~\ref{lem:generic} shows that $W(F)$ is 
no oval for $d=3$, $n\geq 6$. Among real matrices, ovals are generic for 
$d\geq 2$ and $n\leq 2$, but do not appear for $d=3$ and $n\geq 4$. 
We also point out in Section~\ref{sec:boundary} that the discriminant 
vanishes at normal vectors of large faces. This gives an easy to check
condition for the (non-) existence of large faces, because a sum of 
squares decompositions of the modulus of the discriminant 
\cite{Ilyushechkin1992} can be used.\\
%
%
{\par\noindent\footnotesize
{\em Acknowledgements.}
SW thanks Didier Henrion, Konrad Schm\"udgen, and Rainer Sinn for discussions.
It is a pleasure to thank the "Complexity Garage" at the Jagiellonian University,
where all the 3D printouts were made, and to Lia Pugliese for taking their 
photos. We acknowledge support by a Brazilian Capes scholarship (SW), by the 
Polish National Science Center under the project number DEC-2011/02/A/ST1/00119 
(K{\.Z}) and by the project \#56033 financed by the Templeton Foundation.
Research was partially completed while SW was visiting the Institute for 
Mathematical Sciences, National University of Singapore in 2016. }
%
%
\section{Quantum states}
\label{sec:qm}
\par
Our interest in the joint numerical range is its role in quantum 
mechanics where the hermitian matrices $\her_d$ are called Hamiltonians 
or observables, see e.g.\ \cite{BengtssonZyczkowski2006}, and they 
correspond to physical systems with $d$ energy levels or measurable 
quantities having $d$ possible outcomes. 
\par
Usually a (complex) C*-subalgebra $\cA$ of $M_d$ is considered as the 
algebra of observables of a quantum system \cite{AlfsenShultz2001}. If 
$a\in M_d$ is positive semi-definite then we write $a\succeq 0$. The 
physical states of the quantum system are described by $d$-by-$d$ density 
matrices which form the {\em state space} of $\cA$,
\begin{equation}\label{eq:state-space}
\cM(\cA):=\{\rho\in\cA\mid\rho\succeq 0,\tr(\rho)=1\}.
\end{equation}
It is well-known that $\cM(\cA)$ is a compact convex subset of $\her_d$,
see for example Theorem 4.6 of \cite{AlfsenShultz2001}. We are mainly 
interested in $\cM_d:=\cM(M_d)$, but in Sec.~\ref{sec:boundary} also in 
the compressed algebra $pM_dp$ where $p\in M_d$ is a projection, that is 
$p^2=p^*=p$. The state space $\cM(pM_dp)$ is, as we recall in 
Sec.~\ref{sec:boundary}, an exposed face of $\cM_d$, see 
\cite{AlfsenShultz2001,Weis2011}. The state space $\cM_2$ is a Euclidean 
ball, called {\em Bloch ball}, but $\cM_d$ is not a ball 
\cite{BengtssonZyczkowski2006} for $d\geq 3$. Although several attempts were 
made to analyze properties of this set
\cite{JS01,Bengtsson-etal2013,SB13,KKLM16,GSS16}, its complicated structure 
requires further studies.
\par
We use the inner product $\langle a,b\rangle:=\tr(a^*b)$,
$a,b\in M_d$. For any state $\rho\in\cM_d$ and Hamiltonian $a\in\her_d$, the
real number $\langle\rho,a\rangle$ is the expectation value of possible 
outcomes of measurements of $a$. The state $\rho$ is a {\em pure state} if 
$\rho$ is a rank-one projection. The pure state which is the projection onto 
the span of a unit vector $x\in\bC^d$ is denoted by $\rho=|x\rangle\langle x|$ 
and 
\[
\langle\rho,a\rangle 
=
\langle|x\rangle\langle x|,a\rangle 
=\langle x,ax\rangle.
\]
Therefore, the standard numerical range $W(a)$ of a {\em hermitian} matrix $a$ 
is the set of expectation values of $a$ obtained from all pure states. An 
arbitrary state of $\cM_d$, which may not be pure, is called a 
{\em mixed state}. The spectral theorem applied to a mixed state shows 
that the convex hull of $W(a)$ is the set of expectation values of $a$ obtained
from all mixed states. Since $W(a)$ is convex, no convex hull operation is needed 
and therefore $W(a)$ can be identified as a projection of $\cM_d$ onto a line. 
\par
Similarly, the standard numerical range $W(F_1,F_2)$ of a {\em non-hermitian}
operator $F_1+\ii F_2$ is convex. So $W(F_1,F_2)$ is the set of the expectation 
values of measurements of two hermitian operators $F_1$ and $F_2$ performed on 
two copies of the same mixed quantum state. In other words, $W(F_1,F_2)$ is a 
projection of $\cM_d$ onto a two-plane 
\cite{DGHMPZ11,PMGHZ15}. The Dvoretzky theorem \cite{Dv61} implies that for 
large dimension $d$ a generic 2D projection of the convex set $\cM_d$ is 
close to a circular disk, so that the numerical range of a non-hermitian 
random matrix of the Ginibre ensemble typically forms a disk \cite{CGLZ14}.
\par
In this work we analyze joint numerical ranges of triples of hermitian 
matrices of size three. These joint numerical ranges are convex and can be 
interpreted as sets of expectation values of three hermitian observables 
performed on three copies of the same mixed quantum state. They form 
projections of the 8D set of density matrices of size three into a 
three-plane \cite{GutkinZyczkowski2013}. 
\par
An example of projection into high-dimensional planes is the map from the 
states of a composite system to marginals of subsystems. The geometry of 
three-dimensional projections of two-party marginals was recently studied 
\cite{ZDLVHV14,CJLSZ16,Chen-etal2017} to investigate many-body quantum 
systems.
\par
To formalize the discussion of expectation values and projections of 
the set $\cM_d$ we consider arbitrary $d,n\in\bN$ and $n$ hermitian matrices 
$F=(F_1,\ldots,F_n)\in(\her_d)^n$ of size $d$. We will use the linear map 
\[\textstyle
\bE_F:\her_d\to\bR^n,\quad
a\mapsto (\langle a,F_1\rangle,\ldots,\langle a,F_n\rangle)
\]
to study the image 
\begin{equation}\label{eq:defL}\textstyle
L(F):=\bE_F(\cM_d)=\{\bE_F(\rho)\mid\rho\in\cM_d\}\subset\bR^n
\end{equation}
of the state space $\cM_d=\cM(M_d)$ defined in (\ref{eq:state-space}).
The set $L(F)$ was called {\em joint algebraic numerical range} 
\cite{Mueller2010}, also {\em convex support} \cite{Weis2011} in analogy 
with statistics \cite{Barndorff-Nielsen1978}. The compact convex set 
$L(F)$ is the convex hull of the joint numerical range 
\begin{equation}\label{eq:L=convF}\textstyle
L(F)=\conv(W(F)).
\end{equation}
Proofs of equation (\ref{eq:L=convF}) can be found in 
\cite{Mueller2010,GutkinZyczkowski2013}. We recall that $L(F)=W(F)$ holds 
for $n=3$ and $d\geq 3$ where $W(F)$ is convex. In what follows, we will 
work mostly with $L(F)$ rather than $W(F)$.
\par
Some of the 3D images shown in Figure~\ref{fig:print} are generated using 
random sampling -- this method is simple conceptually and produces objects 
which are accurate enough for use in printing. In this numerical procedure, 
which we implemented in Mathematica, we calculate a finite number (of order 
of $10^5$) of points inside $W(F)$,
\[
\{(\langle x,F_1x\rangle,\ldots,\langle x,F_nx\rangle)
\mid x\in S\},
\]
where $S$ is the set of points sampled from the uniform distribution on 
the unit sphere of $\bC^d$ (this step is realized by sampling 
points from complex $d$-dimensional Gaussian distribution and normalizing 
the result). A convex hull of generated points is then calculated using 
\verb|ConvexHullMesh| procedure and exported to an \verb|.stl| file, which 
contains a description of the 3D object recognized by the software used 
in printing. The final objects were made with PIRX One 3D printer.
%
%
\section{Normal cones and ovals}
\label{sec:lattices}
\par
We show that joint algebraic numerical ranges have in a sense many 
normal cones. We prove that this property allows to characterize ovals 
in terms of strict convexity.
\par
A {\em face} of a convex subset $C\subset\bR^m$, $m\in\bN$, is a convex 
subset of $C$ which contains the endpoints of every open segment in 
$C$ which it intersects. An {\em exposed face} of $C$ is defined as 
the set of maximizers of a linear functional on $C$. If $C$ is non-empty 
and compact, then for every $u\in\bR^m$ the set 
\begin{equation}\label{eq:expface}
\bF_C(u):=\argmax_{x\in C}\langle x,u\rangle
\end{equation}
is an exposed face of $C$. By definition, the empty set is also an 
exposed face (then the set of exposed faces forms a lattice).
It is well-known that every exposed face is a face. If a face 
(resp.~exposed face) is singleton, then we call its element an 
{\em extreme point} (resp.~{\em exposed point}). A face 
(resp.~exposed face) of $C$ which is different from $\emptyset,C$ 
is called {\em proper face} (resp.~proper exposed face).
\par
Let $C\subset\bR^m$ be a convex subset and $x\in C$. The {\em normal cone} 
of $C$ at $x$ is 
\[
N(x):=\{u\in\bR^m\mid\forall y\in C:\,\langle u,y-x\rangle\leq 0\}.
\] 
Elements of $N(x)$ are called {\em (outer) normal vectors} of $C$ at $x$. 
It is well-known that there is a non-zero normal vector of $C$ at $x$ if 
and only if $x$ is a boundary point of $C$. In that case $x$ is 
{\em smooth} if $C$ admits a unique outer unit normal vector at $x$. 
\par
The normal cone of $C$ at a non-empty face $G$ of $C$ is well-defined as 
the normal cone $N(G):=N(x)$ of $C$ at any point $x$ in the relative 
interior of $G$ (the relative interior of $G$ is the interior of $G$ with 
respect to the topology of the affine hull of $G$). See for example 
Section~4 of \cite{Weis2012a} about the consistency of this definition,
and set $N(\emptyset):=\bR^m$. The convex set $C$ is a 
{\em convex cone} if $C\neq\emptyset$ and if $x\in C$, $\lambda\geq 0$ 
implies $\lambda x\in C$. A {\em ray} is a set of the form
$\{\lambda\cdot u\mid \lambda\geq 0\}\subset\bR^m$ for non-zero 
$u\in\bR^m$. An {\em extreme ray} of $C$ is a ray which is a face of $C$.
\par
We denote the set of exposed faces and normal cones of $C$ by $\cE_C$ and 
$\cN_C$, respectively. Each of these sets is partially ordered by inclusion
and forms a {\em lattice}, that is the infimum and supremum of each pair of 
elements exist. A {\em chain} in a lattice is a totally ordered subset, the 
{\em length} of a chain is the cardinality minus one. The {\em length} of a 
lattice is the supremum of the lengths of all its chains. Lattices of faces 
have been studied earlier \cite{Barker1973,LoewyTam1986}, in particular 
these of state spaces \cite{AlfsenShultz2001}, and linear images $L(F)$ of 
state spaces \cite{Weis2011}. By Proposition~4.7 of \cite{Weis2012a}, if $C$ 
is not a singleton then
\begin{equation}\label{eq:exp-norm-iso}
\cE_C\to\cN_C,
\quad
G\mapsto N(G)
\end{equation}
is an antitone lattice isomorphism. This means that the map is a bijection 
and for all exposed faces $G,H$ we have $G\subset H$ if and only if 
$N(G)\supset N(H)$.
\par
What makes a joint algebraic numerical range $C$ special is that all non-empty 
faces of its normal cones are normal cones of $C$, too, as we will see in 
Lemma~\ref{lem:N=T}. For two-dimensional $C\subset\bR^2$ this means that a 
boundary point of $C$ is smooth unless it is the intersection of two 
one-dimensional faces of $C$, as one can see from the isomorphism 
(\ref{eq:exp-norm-iso}). That property is well-known \cite{Bebiano1986} for the 
numerical range $W(F_1,F_2)$ of a matrix $A=F_1+\ii F_2\in M_d$. For example, 
the half-moon $\{z\in\bC : |z|\leq 1, \Re(z)\geq 0\}$ is not the numerical 
range of any matrix. This also follows from Anderson's theorem 
\cite{ChienNakazato1999} which asserts that if $W(F_1,F_2)$ is included in the 
unit disk and contains $d+1$ distinct points of the unit circle, then 
$W(F_1,F_2)$ is the unit disk.
\par
To prove the lemma we introduce the Definitions~6.1 and~7.1 of 
\cite{Weis2012a} for the special case of a non-empty, compact, and convex 
subset $C\subset\bR^m$. Let $u\in\bR^m$ be a non-zero vector. Then $u$ is 
called {\em sharp normal} for $C$ if for every relative interior point $x$ 
of the exposed face $\bF_C(u)$ the vector $u$ is a relative interior point 
of the normal cone of $C$ at $x$. The {\em touching cone} of $C$ at $u$ is 
defined to be the face of the normal cone of $C$ at $\bF_C(u)$ which 
contains $u$ in its relative interior \cite{Schneider2014}. The linear 
space $\bR^m$ and the orthogonal complement of the translation vector space 
of the affine hull of $C$ are touching cones of $C$ by definition. We point 
out that every normal cone of $C$ is a touching cone of $C$.
\begin{lem}\label{lem:N=T}
 Every non-empty face of every normal cone of $L(F)$ is a 
normal cone of $L(F)$. 
\end{lem}
\par\noindent
{\em Proof:}
Propositions 2.9 and 2.11 of \cite{Weis2011} prove that every non-zero 
hermitian $d$-by-$d$ matrix is sharp normal for the state space $\cM_d$. 
Therefore, Proposition~7.6 of \cite{Weis2012a} shows that every touching 
cone of $\cM_d$ is a normal cone of $\cM_d$. Corollary 7.7 of 
\cite{Weis2012a} proves that $L(F)$, being a projection of $\cM_d$, has 
the analogous property that every touching cone of $L(F)$ is a normal 
cone of $L(F)$. The characterization of touching cones as the non-empty 
faces of normal cones, given in Theorem~7.4 of \cite{Weis2012a}, 
completes the proof.
\hspace*{\fill}$\square$\\
\par
We define an {\em oval} as a convex and compact subset of $\bR^m$ with 
interior points each of whose boundary points is a smooth exposed point. 
Notice that ovals are strictly convex. For the following class of convex 
sets strict convexity implies smoothness.
\begin{lem}\label{lem:char-ovals}
Let $C\subset\bR^m$ be a convex and compact subset of $\bR^m$ with 
interior points, such that every extreme ray of every normal cone of $C$ 
is a normal cone of $C$. Then $C$ is an oval if and only if all proper 
exposed faces of $C$ are singletons.
\end{lem}
{\em Proof:}
We assume first that $C$ is an oval. By definition, the boundary of $C$ 
is covered by extreme points. Since $C$ is the disjoint union of the 
relative interiors of its faces, see for example 
Theorem 2.1.2 of \cite{Schneider2014}, this shows that all proper faces of 
$C$ are singletons.
\par
Conversely, we assume that all proper exposed faces of $C$ are singletons.
Since $C$ has full dimension, the proper faces of $C$ cover the boundary 
$\partial C$. Every proper face lies in a proper exposed face, see for 
example Lemma 4.6 of \cite{Weis2012a}, so $\partial C$ is covered by exposed 
points. Let $x$ be an arbitrary exposed point of $C$. We have to show that 
$x$ is a smooth point. As $\dim(C)=m$, the normal cone $N(x)$ contains no 
line and so it has at least one extreme ray which we denoted by $r$ (see 
e.g.~Theorem 1.4.3 of \cite{Schneider2014}). By assumption, $r$ is a normal 
cone of $C$. So 
\[
\{0\}\subset r \subset N(x) \subset \bR^m
\]
is a chain in the lattice $\cN_C$ of normal cones. Thereby the inclusion
$\{0\}\subset r$ is proper. By the antitone isomorphism 
(\ref{eq:exp-norm-iso}), there is an exposed face $F$ with normal cone $r$, 
\[
C \supset F \supset \{x\} \supset \emptyset
\]
is a chain in the lattice $\cE_C$ of exposed faces, and the inclusion
$C \supset F$ is proper. By assumption, all proper exposed faces of $C$ 
are singletons. So $F=\{x\}$ follows. Using the isomorphism 
(\ref{eq:exp-norm-iso}) a second time gives $r=N(x)$, that is 
$x$ is a smooth point.
\hspace*{\fill}$\square$\\
%
%
%
\section{Exposed faces}
\label{sec:boundary}
\par
This section collects methods to study exposed faces of the 
joint algebraic numerical range $L(F)$. We start with the well-known 
representation of exposed faces in terms of eigenspaces of the greatest 
eigenvalues of real linear combinations of $F_1,\ldots,F_n$. This allows 
us to show that the generic shape of $L(F)$ is an oval for $n=1,2,3$
($n=1,2$ for real symmetric $F_i$'s). For $3$-by-$3$ matrices we 
discuss the discriminant of the characteristic polynomial and the 
sum of squares decomposition of its modulus. We further discuss 
pre-images of exposed points. This allows us to prove that $L(F)$ 
is no oval for $d=3$ if $n\geq 6$ ($n\geq 4$ for real symmetric 
matrices). Finally we address corner points.
\par
Let $d,n\in\bN$ be arbitrary. As before we write 
$F=(F_1,\ldots,F_n)\in(\her_d)^n$ and we define
\[
F(u):=u_1F_1+\cdots+u_nF_n,
\qquad
u\in\bR^n.
\] 
By (\ref{eq:defL}) the joint algebraic numerical range $L(F)$ is the image 
of the state space $\cM_d$ under the map $\bE_F$. So all subsets of $L(F)$ 
are equivalently described in terms of their pre-images under the restricted 
map $\bE_F|_{\cM_d}$. In particular, the exposed face $\bF_{L(F)}(u)$ of 
$L(F)$, in the notation from (\ref{eq:expface}), has the pre-image
\begin{equation}\label{eq:face-lift}\textstyle
\bE_F|_{\cM_d}^{-1}(\bF_{L(F)}(u))
=\bF_{\cM_d}(F(u))
=\argmax_{\rho\in\cM_d}\langle\rho,F(u)\rangle.
\end{equation}
See for example Lemma~5.4 of \cite{Weis2012a} for this simple observation. 
The equation (\ref{eq:face-lift}) offers an algebraic description of exposed 
faces of $L(F)$. For $a\in\her_d$ the exposed face $\bF_{\cM_d}(a)=\cM(pM_dp)$ 
of the state space $\cM_d=\cM(M_d)$ is the state space of the 
algebra $pM_dp$ where $p$ is the spectral projection of $a$ 
corresponding to the greatest eigenvalue, see \cite{AlfsenShultz2001} 
or \cite{Weis2011}. Therefore (\ref{eq:face-lift}) shows 
\begin{equation}\label{eq:pre-image}
\bE_F|_{\cM_d}^{-1}(\bF_{L(F)}(u))
=\cM(pM_dp),
\qquad u\in\bR^n,
\end{equation}
where $p$ is the spectral projection of $F(u)$ corresponding to the 
greatest eigenvalue.
\begin{rem}[Spectral representation of faces]\label{rem:h-L}
A proof is given in Section~3.2 of \cite{Gutkin-etal2004} that for 
$u\in\bR^n$ the {\em support function} 
$h_{W(F)}:=\max\{\langle x,u\rangle\mid x\in W(F)\}$ of $W(F)$ is the 
greatest eigenvalue of $F(u)$. This result goes back to Toeplitz 
\cite{Toeplitz1918} for $n=2$. The same conclusion follows also from 
(\ref{eq:L=convF}) and (\ref{eq:pre-image}), in particular 
$h_{W(F)}(u)=\max\{\langle x,u\rangle\mid x\in L(F)\}$.
\end{rem}
\par
The generic joint algebraic numerical range of at most three hermitian 
matrices is an oval. 
\begin{thm}\label{thm:generic}
Let $n\in\{1,2,3\}$ and $d\geq2$. Then the set of $n$-tuples of hermitian 
$d$-by-$d$ matrices $F\in(\her_d)^n$ such that $L(F)$ is an oval is open 
and dense in $(\her_d)^n$.
\end{thm}
{\em Proof:}
For $n=1,2,3$ and $d\in\bN$ the set $\mathcal{O}_1$ of all $F\in(\her_d)^n$ 
where every matrix in the pencil $\{F(u)\mid u\in\bR^n\setminus\{0\}\}$
has $d$ simple eigenvalues is open and dense in $(\her_d)^n$, this was shown 
in Prop.~4.9 of \cite{Gutkin-etal2004}. Hence, for $F\in\mathcal{O}_1$ all 
proper exposed faces of $L(F)$ are singletons by (\ref{eq:pre-image}). 
Secondly, since $n+1\leq\dim_\bR(\her_d)=d^2$ holds by the assumptions 
$n\leq 3$ and $d\geq 2$, it is easy to prove that $\id_d,F_1,\ldots,F_n$ are 
linearly independent for $F$ in an open and dense subset $\mathcal{O}_2$ of 
$(\her_d)^n$, that is $\dim(L(F))=n$ holds for $F\in\mathcal{O}_2$. The 
extreme rays of every normal cone of $L(F)$ are normal cones of $L(F)$ by 
Lemma~\ref{lem:N=T}. Hence Lemma~\ref{lem:char-ovals} proves that $L(F)$ is 
an oval for all $F$ in $\mathcal{O}_1\cap\mathcal{O}_2$. The proof is 
completed by observing that the intersection of two open and dense subsets 
of any topological space is open and dense.
\hspace*{\fill}$\square$\\
\par
Let us now focus on $3$-by-$3$ matrices ($d=3$). As explained earlier in this
section, every proper exposed face of the state space $\cM_3=\cM(M_3)$ is the 
state space $\cM(pM_3p)$ of the algebra $pM_3p$ for a projection $p\in M_3$ of 
rank one or two. In the former case $\cM(pM_3p)$ is a singleton and in the 
latter case a three-dimensional Euclidean ball. Hence (\ref{eq:face-lift}) 
shows that every proper exposed face of $L(F)$ is a singleton, segment, filled 
ellipse, or filled ellipsoid.
\begin{lem}\label{lem:large-faces-exposed}
Let $F$ be an $n$-tuple of hermitian $3$-by-$3$ matrices. Then every proper 
face of $L(F)$ is a singleton, segment, filled ellipse, or filled ellipsoid. 
If that face is no singleton then it is an exposed face of $L(F)$.
\end{lem}
{\em Proof:}
Every proper face $G$ of $L(F)$ lies in a proper exposed face $H$ of $L(F)$ 
(see for example Lemma 4.6 of~\cite{Weis2012a}), hence $G\subset H$ is a 
face of $H$. As mentioned above, $H$ is a singleton, segment, ellipse, or 
ellipsoid. Therefore $G=H$ holds if $G$ is no singleton.
\hspace*{\fill}$\square$\\
\par
The next aim is to provide a method to certify that all large faces,
defined in (\ref{eq:large-face}) for $d=n=3$, were found. To this end
we use the discriminant and a sum of squares decomposition of its 
modulus. 
\begin{rem}[Discriminant method]\label{rem:discrim-method}
Recall from (\ref{eq:defL}) that $L(F)=\bE_F(\cM_3)$ is a projection of a 
state space. Hence, if the exposed face $\bF_{L(F)}(u)$ of $L(F)$, defined 
by $u\in\bR^3$, is a large face, then its pre-image 
$\bE_F|_{\cM_d}^{-1}(\bF_{L(F)}(u))$ is necessarily no singleton. As we 
pointed out in (\ref{eq:pre-image}) this means that the greatest eigenvalue 
of $F(u)$ is degenerate, which is equivalent to a vanishing discriminant as 
we see next.
\end{rem}
\par
Let $a_1,a_2,a_3\in\bC$ and consider the polynomial 
$p(\lambda)=-\lambda^3+a_1\lambda^2+a_2\lambda+a_3$ of degree three. 
The {\em discriminant} of $p$, see Section A.1.2 of \cite{Fischer2001}, is
\[
-(27 a_3^2 + 18 a_1 a_2 a_3 - 4 a_1^3 a_3 + 4 a_2^3 - a_1^2 a_2^2).
\]
Let $\lambda_1,\lambda_2,\lambda_3\in\bC$ denote the roots of $p$. Then 
the discriminant of $p$ can be written
\[\textstyle
\Pi_{1\leq i<j\leq 3}(\lambda_i-\lambda_j)^2.
\]
The {\em discriminant} $\delta(A)$ of a 3-by-3 matrix $A\in M_3$ is the 
discriminant of the characteristic polynomial $\det(A-\lambda\id)$. 
So, $A$ has a multiple eigenvalue if and only if $\delta(A)=0$. 
\par
Let $Z\in M_3$ be a normal $3$-by-$3$ matrix, that is $Z^*Z=ZZ^*$. The 
entries of the matrices $Z^0=\id$, $Z^1=Z$, and $Z^2=Z Z$ can be combined 
into a $9$-by-$3$ matrix $Z_*$ by choosing an ordering of 
$\{1,2,3\}^{\times2}=\{1,2,3\}\times\{1,2,3\}$. The $i$-th column of $Z_*$ 
is defined to be equal to $Z^i$ in that ordering for $i=0,1,2$. Now the 
absolute value of the discriminant $\delta(Z)$ is \cite{Ilyushechkin1992} 
\begin{equation}\label{eq:Ilyushechkin}\textstyle
|\delta(Z)|=\sum_\nu |M_\nu|^2,
\end{equation}
where the sum extends over the $84$ subsets $\nu\subset\{1,2,3\}^{\times2}$ 
of cardinality three and where $M_\nu$ is the $3$-by-$3$ minor of the rows 
of $Z_*$ which are indexed by $\nu$. The theory of discriminants or 
(\ref{eq:Ilyushechkin}) show that $|\delta(Z)|$ is homogeneous of degree 
six. It is worth noting that the discriminant of a real symmetric $3$-by-$3$ 
matrix can be decomposed into a sum of five squares \cite{Domokos2011}.
\par
For $3$-by-$3$ matrices, a vanishing discriminant is not only necessary 
(Remark~\ref{rem:discrim-method}) but also almost sufficient for the existence
of large faces, with the exception of special Euclidean balls. To describe 
this problem more precisely, let $d,n\in\bN$ be arbitrary. We define an 
equivalence relation on $(\her_d)^n$. For $F=(F_1,\ldots,F_n)\in(\her_d)^n$ 
and a unitary $U\in M_d$ let $U^*FU:=(U^*F_1U,\ldots,U^*F_nU)$. Two tuples 
$F,G\in(\her_d)^n$ are equivalent if and only if either
\begin{equation}\label{eq:equiv-cond1}\textstyle
\mbox{$G=U^*FU$ holds for some unitary $U\in M_d$}
\end{equation}
or
\begin{equation}\label{eq:equiv-cond2}\textstyle
\mbox{$\id_d,G_1,\ldots,G_n$ and $\id_d,F_1,\ldots,F_n$ have the same span.}
\end{equation}
The statement (\ref{eq:equiv-cond1}) means that the equivalence classes are 
invariant under unitary similarity with any unitary $U\in M_d$,
\begin{equation}\label{eq:trafo1}\textstyle
(\her_d)^n\to(\her_d)^n,
\quad
F\mapsto U^*FU.
\end{equation}
Joint algebraic numerical ranges are fixed under these maps, $L(U^*FU)=L(F)$. 
The statement (\ref{eq:equiv-cond2}) means that the equivalence classes are 
invariant under the action of the affine group of $\bR^n$. More precisely, 
let $A=(a_{i,j})\in\bR^{n\times n}$ be an invertible matrix and $b\in\bR^n$. 
There are two affine transformations
\begin{equation}\label{eq:trafo2a}\textstyle
\alpha:\bR^n\to\bR^n, 
\quad
x\mapsto(\sum_{j=1}^n a_{i,j}x_j+b_i)_{i=1}^n
\end{equation}
and 
\begin{equation}\label{eq:trafo2}\textstyle
\beta:(\her_d)^n\to(\her_d)^n,
\quad 
F\mapsto(\sum_{j=1}^n a_{i,j}F_j+b_i\id_d)_{i=1}^n.
\end{equation}
Notice that $\alpha(L(F))=L(\beta(F))$ holds, as 
$\alpha\circ\bE_F(X)=\bE_{\beta(F)}(X)$ for all $X\in\her_d$ of trace one. 
In other words, joint algebraic numerical ranges and $n$-tuples of 
hermitian matrices transform equivariantly under the affine group of 
$\bR^n$.
\par
We call $(v_1,\ldots,v_k)$ a {\em real $k$-frame} of $\bC^2$ if 
$v_1,\ldots,v_k\in\bC^2$ are real linearly independent. The tuple 
$(v_1,\ldots,v_k)$ is an orthonormal real $k$-frame of $\bC^2$ if 
$v_1,\ldots,v_k\in\bC^2$ are orthonormal with respect to the Euclidean 
scalar product which is the real part of the standard inner product of 
$\bC^2$.
\begin{lem}\label{lem:M2-singletons}
Let $n\in\bN$, $F\in(\her_3)^n$, $D:=\dim(L(F))$, $D\geq 1$, and assume 
that the pre-image of some exposed point of $L(F)$ under $\bE_F|_{\cM_3}$ 
is no singleton. 
\par
Then $D\leq 5$ and $F$ is equivalent modulo (\ref{eq:trafo1}) and 
(\ref{eq:trafo2}) to $G=(G_1,\ldots,G_n)$ where
\begin{equation}\label{eq:Gdef}
G_1=\left(\begin{array}{ccc}
1 & 0 & 0\\
0 & 1 & 0\\
0 & 0 & -1
\end{array}\right),
\qquad
G_i=\left(\begin{array}{cc|c}
0 & 0 & \raisebox{-7pt}{$v_i$}\\[-5pt]
0 & 0 & \\ \hline
\multicolumn{2}{c|}{v_i^*} & 0
\end{array}\right),
\quad 2\leq i\leq D,
\end{equation}
for a real $(D-1)$-frame $(v_2,\ldots,v_D)$ of $\bC^2$, and where $G_i=0$ 
for $D<i\leq n$. The real frame may be chosen to be orthonormal. 
Specifically, if $D\leq 4$ then there are $\varphi\in[0,\pi]$ and 
$\theta\in[0,2\pi)$ such that $v_2,\ldots,v_D$ can be taken from the list
\begin{equation}\label{eq:vectorsD4}
v_2=
\left(\begin{smallmatrix}
1 \\
0
\end{smallmatrix}\right),
\quad
v_3=
\left(\begin{smallmatrix}
\ii\cos(\varphi) \\
\sin(\varphi)
\end{smallmatrix}\right),
\quad\mbox{and}\quad
v_4=
\cos(\theta)
\left(\begin{smallmatrix}
-\ii\sin(\varphi) \\
\cos(\varphi)
\end{smallmatrix}\right)
+\sin(\theta)
\left(\begin{smallmatrix}
0 \\
\ii 
\end{smallmatrix}\right).
\end{equation}
Thereby $\varphi$ is unique if $D=3$. Both $\varphi$ and $\theta$ are 
unique if $D=4$. If $D=5$ then one can take 
\begin{equation}\label{eq:vectorsD5}
v_2=
\left(\begin{smallmatrix}
1 \\
0
\end{smallmatrix}\right),
\quad
v_3=
\left(\begin{smallmatrix}
\ii\\
0
\end{smallmatrix}\right),
\quad
v_4=
\left(\begin{smallmatrix}
0 \\
1
\end{smallmatrix}\right),
\quad\mbox{and}\quad
v_5=
\left(\begin{smallmatrix}
0 \\
\ii
\end{smallmatrix}\right).
\end{equation}
If the matrices $F$ are real symmetric, then $D\leq 3$ and one can 
choose $v_2,\ldots,v_D$ from the list
\begin{equation}\label{eq:vectorsD3real}
v_2=
\left(\begin{smallmatrix}
1 \\
0
\end{smallmatrix}\right)
\quad\mbox{and}\quad
v_3=
\left(\begin{smallmatrix}
0 \\
1
\end{smallmatrix}\right).
\end{equation}
For all matrix tuples $G$ of the form (\ref{eq:Gdef}) with $(D-1)$-frames 
(\ref{eq:vectorsD4}), (\ref{eq:vectorsD5}), or (\ref{eq:vectorsD3real}), 
the joint algebraic numerical range $L(G)$ is the cartesian product of the 
unit ball $B^D=\{y\in\bR^D\mid y_1^2+\cdots+y_D^2\leq 1\}$ of $\bR^D$ with 
the origin of $\bR^{n-D}$. The pre-image $\bE_G|_{\cM_3}^{-1}(1,0,\ldots,0)$ 
is a three-dimensional Euclidean ball and the pre-images of all other exposed 
points of $L(G)$ are singletons. 
\end{lem}
{\em Proof:}
Let $x$ be an exposed point of $L(F)$ with multiple pre-images, say 
$u\in\bR^n$ is a unit vector and $\{x\}=\bF_{L(F)}(u)$. Applying a rotation 
(\ref{eq:trafo2}) of $\bR^n$ we take $u:=(1,0,\ldots,0)$. By 
(\ref{eq:pre-image}) the pre-image of $x$ is 
$\bE_F|_{\cM_3}^{-1}(x)=\cM(p\cM_3p)$ where $p$ is the spectral projection 
of $F(u)=F_1$ corresponding to the greatest eigenvalue of $F_1$. Since 
\[
\bE_F|_{\cM_3}^{-1}(x)\neq\{\tfrac{p}{\tr(p)}\},\cM_3,
\] 
it follows that $p$ has rank two. Notice that $pF_1p,\ldots,pF_np$ are scalar 
multiples of $p$. Otherwise there will be $\rho_1,\rho_2\in\cM(p\cM_3p)$ and 
an index $i\in\{1,\ldots,n\}$ such that 
$\langle\rho_1,F_i\rangle\neq\langle\rho_2,F_i\rangle$, but this contradicts 
the assumption that $\{x\}=\bE_F(\cM(p M_3p))$ is a singleton. A unitary 
similarity (\ref{eq:trafo1}) and another affine map (\ref{eq:trafo2}) 
transform $F$ into the tuple $G$ defined in (\ref{eq:Gdef}).
\par
The real $(D-1)$-frame $\chi:=(v_2,\ldots,v_D)$ may be transformed into an
orthogonal real frame using the unitary group $U(2)$ and the general linear 
group ${\rm GL_4}(\bR)$ acting on $\bC^2\cong\bR^4$. More precisely, a 
unitary $V\in U(3)$ which acts on $G$ {\em via} (\ref{eq:trafo1}) and which 
keeps $G_1$ fixed has the form
\[
V=\left(\begin{smallmatrix}
U & 0 \\
0 & e^{-\ii\phi}
\end{smallmatrix}\right)
\]
for some unitary $U\in U(2)$ and $\phi\in\bR$. The action of $V$ on
$G_2,\ldots,G_D$ is 
\[
V 
\left(\begin{smallmatrix}
0 & v_i \\
v_i^* & 0 
\end{smallmatrix}\right)
V^* 
= 
\left(\begin{smallmatrix}
0 & e^{\ii\phi}Uv_i \\
(e^{\ii\phi}Uv_i)^* & 0 
\end{smallmatrix}\right).
\] 
An affine map (\ref{eq:trafo2}) which fixes 
$G_1$ and $G_{D+1}=\cdots=G_n=0$ acts on $\chi$ by taking invertible real 
linear combinations. So the general linear group ${\rm GL_4}(\bR)$ acts on 
$\chi$. If $F$ is real symmetric, then the orthogonal group 
$O(2)\subset U(2)$ suffices. These group actions lead to the orthogonal 
real frames (\ref{eq:vectorsD4}), (\ref{eq:vectorsD5}), and 
(\ref{eq:vectorsD3real}).
\par
Let us analyze $L(G)$. Since $G_i=0$ for $i>D$, it suffices to study $D=n$.
Remark~\ref{rem:h-L} shows that for $u\in\bR^n$ 
\[
h(u):=\max_{x\in L(G)}\langle u,x\rangle
\] 
is the greatest eigenvalue of $G(u)=u_1G_1+\cdots+u_nG_n$. An easy 
computation shows that if $u$ is a unit vector then the matrix $G(u)$ has 
eigenvalues $\{-1,u_1,1\}$. So $h(u)=1$ holds for all unit vectors 
$u\in\bR^n$ and this shows that $L(G)$ is the unit ball $B^n$. Let 
$v:=(1,0,\ldots,0)\in\bR^n$. The pre-image $\bE_{G}|_{\cM_3}^{-1}(v)$ of 
the exposed point $\{v\}=\bF_{L(G)}(v)$ of $L(G)$ is a three-dimensional 
ball since the greatest eigenvalue of $G_1$ is degenerate 
(\ref{eq:pre-image}). As we point out in Remark~\ref{rem:discrim-method},
to see that $v$ is the unique exposed point of $L(G)$ with multiple 
pre-image points, it suffices to show that the discriminant 
$\delta(G(u))$ is non-zero for unit vectors $u\in\bR^n$ which are not 
collinear with $v$. But this follows from the formula 
\[
\delta(G(u))=4(u_2^2 + \cdots + u_n^2)^2
\]
which is readily verified.
\hspace*{\fill}$\square$\\
\par
It is worth to remark on unitary (ir-) reducibility in the context
of pre-images.
\begin{rem}
Let $n\in\bN$, $F\in(\her_3)^n$, and let $L(F)$ have an exposed point with 
multiple pre-images under $\bE_F|_{\cM_3}:\cM_3\to L(F)$. It was shown in 
Theorem~3.2 of \cite{Leake-etal2014} for such tuples $F$ that if the dimension
$D=\dim(L(F))$ is at most $D=2$ then $F$ is unitarily reducible. The same 
conclusion can be drawn from Lemma~\ref{lem:M2-singletons}. With rare 
exceptions, the lemma shows also that if $D\geq 3$ then $F$ is unitarily 
irreducible. The exceptions are those $F$ where $D=3$ and where $F$ is 
equivalent modulo (\ref{eq:trafo1}) and (\ref{eq:trafo2}) to an $n$-tuple 
$G$ with vectors 
\[
v_2=
\left(\begin{smallmatrix}
1 \\
0
\end{smallmatrix}\right)
\quad\mbox{and}\quad
v_3=
\left(\begin{smallmatrix}
\pm\ii\\
0
\end{smallmatrix}\right)
\]
specified in equation (\ref{eq:vectorsD4}).
\end{rem}
\par
While $L(F)$ is generically an oval for $d\geq 2$ and $n\leq 3$ (by 
Theorem~\ref{thm:generic}), we now exclude ovals for $d=3$ and large 
$n$.
\begin{lem}\label{lem:generic}
Let $F\in(\her_3)^n$ and $D=\dim(L(F))$. If $D\geq 6$ ($D\geq 4$ 
suffices if the $F_i$'s are real symmetric), then $L(F)$ is no oval.
\end{lem}
{\em Proof:}
Since $n\geq D$ holds, the bound on $D$ implies $n\geq 4$ (resp.~$n\geq 3$ 
for real symmetric matrices). Thus, Theorem D (resp.~Theorem B) of 
\cite{Friedland-etal1984} proves that there is a non-zero $u\in\bR^n$ such 
that $F(u)=u_1F_1+\cdots+u_nF_n$ has a multiple eigenvalue. So the greatest 
eigenvalue of $F(v)$ is degenerate, either for $v=u$ or for $v=-u$. As we 
see from (\ref{eq:pre-image}), this means that the exposed face 
$\bF_{L(F)}(v)$ has multiple pre-image points under $\bE_F|_{\cM_3}$. If 
$L(F)$ is an oval, then $\bF_{L(F)}(v)$ is a singleton and then 
Lemma~\ref{lem:M2-singletons} shows $D\leq 5$ ($D\leq 3$ if the matrices $F$ 
are real symmetric).
\hspace*{\fill}$\square$\\
\par
We finish the section with an analysis of corner points of a 
convex compact subset $C\subset\bR^m$. A point $x\in C$ is a 
{\em corner point} \cite{Fiedler1981} of $C$ if the 
normal cone $N(x)$ of $C$ at $x$ has dimension $m$. A point 
$x\in C$ is a {\em conical point} \cite{BindingLi1991} of $C$ if 
$C\subset x+K$ holds for a closed convex cone $K\subset\bR^m$ 
containing no line. The {\em polar} of a closed convex cone 
$K\subset\bR^m$ is
\[
K^\circ=\{u\in\bR^m\mid\forall x\in K:\langle x,u\rangle\leq 0\}.
\]
We recall that $K^\circ$ is a closed convex cone and
$K=(K^\circ)^\circ$, see for example \cite{Schneider2014}.

\begin{lem}\label{eq:conical-corner}
Let $C\subset\bR^m$ be a convex compact subset and $x\in C$. 
Then $x$ is a conical point of $C$ if and only if $x$ is a 
corner point $C$.
\end{lem}
{\em Proof:}
For any point $x\in C$, the smallest closed convex cone containing 
$C-x$ is the polar $N(x)^\circ$ of the normal cone $N(x)$ of
$C$ at $x$, see for example equation (2.2) of \cite{Schneider2014}. 
So for an arbitrary closed convex cone $K\subset\bR^m$ we have 
$C-x\subset K\iff N(x)^\circ\subset K$, that is 
\[
C\subset x+K\iff K^\circ\subset N(x).
\]
The observation that $K$ contains no line if and only if $K^\circ$ 
has full dimension $m$ then proves the claim.
\hspace*{\fill}$\square$\\
\par
The existence of corner points of $L(F)$ has strong algebraic 
consequences for $F$.
\begin{lem}\label{lem:BindingLi}
Let $F\in(\her_3)^n$, and let $p$ be a corner point of $L(F)$. Then 
$F$ is unitarily reducible and there exists a non-zero vector 
$x\in\bC^d$ such that $F_ix=p_ix$ holds for $i=1,\ldots,n$.
\end{lem}
{\em Proof:}
The equivalence of the notions of {\em conical point} and 
{\em corner point} is proved in Lemma~\ref{eq:conical-corner}. The 
remaining claims are proved in Proposition 2.5 of \cite{BindingLi1991}.
\hspace*{\fill}$\square$\\
\par
We derive a classification of corner points of $L(F)$ for 3-by-3 matrices.
\begin{lem}\label{lem:corner-of-3-by-3}
Let $F\in(\her_3)^n$, $D:=\dim(L(F))$, and let $p\in\bR^n$ be a corner 
point of $L(F)$. Then $D\in\{0,1,2,3,4\}$ and, ignoring $D=0,1$, the 
joint algebraic numerical range $L(F)$ is the convex hull of the union 
of $\{p\}$
\begin{itemize}
\item ($D=2$)
with a segment whose affine hull does not contain $p$ 
or with an ellipse which contains $p$ in its affine hull 
but not in its convex hull,
\item ($D=3$)
with an ellipse whose affine hull does not contain $p$ 
or with an ellipsoid which contains $p$ in its affine hull 
but not in its convex hull,
\item ($D=4$)
with an ellipsoid whose affine hull does not contain $p$.
\end{itemize}
\end{lem}
{\em Proof:} 
Lemma~\ref{lem:BindingLi} proves that there exists a unitary $U\in M_3$ 
such that $U^*FU$ has the block-diagonal form 
$U^*FU=\left( \,(p_1)\oplus G_1,\ldots,(p_n)\oplus G_n\,\right)$ with 
$G\in(\her_2)^n$. The joint algebraic numerical range $L(F)=L(U^*FU)$ is 
the convex hull of the union of $L(G)$ and $\{p\}$. Since $L(G)$ is a 
singleton, a segment, a filled ellipse, or a filled ellipsoid, only the 
cases listed above do occur.
\hspace*{\fill}$\square$\\
%
%
%
%
\section{Arguments for the classification}
\label{sec:classification}
\par
Details of Theorem~\ref{thm:intro} are discussed concerning intersections 
of large faces, a graph embedding, and corner points.
\par
We consider the joint numerical range $L(F)$ of a triple 
$F=(F_1,F_2,F_3)\in(\her_3)^3$ of hermitian $3$-by-$3$ matrices. We recall 
from (\ref{eq:large-face}) that a large face of $L(F)$ is a proper exposed 
face of $L(F)$ which is no singleton. Equivalently, a large face is an 
exposed face of $L(F)$ of the form of a segment or of the form of an 
ellipse, but different from $L(F)$ itself. The set of large faces of $L(F)$ 
is denoted by $\cL(F)$. Let 
\begin{equation}\label{eq:projection-x-y}
P:\bR^3\to\bR^2, 
\quad
(x_1,x_2,x_3)\mapsto(x_1,x_2)
\end{equation}
denote the projection onto the $x_1$-$x_2$-plane.
\begin{lem}\label{lem:class1}
Let $F\in(\her_3)^3$, let $G_1,G_2\in\cL(F)$, and let $G_1\neq G_2$. Then 
$G_1$ and $G_2$ intersect in a unique point which is an extreme point of 
$G_1$ and of $G_2$. There is a rotation $\alpha$ of $\bR^3$ in the 
notation of (\ref{eq:trafo2a}) and a corresponding map $\beta$ defined in 
(\ref{eq:trafo2}), such that $\alpha(L(F))=L(\beta(F))$ and such that 
$P(\alpha(G_1))$ and $P(\alpha(G_2))$ are different one-dimensional faces 
of $L(F'_1,F'_2)$ where $(F'_1,F'_2,F'_3)=\beta(F)$. 
\end{lem}
{\em Proof:} 
The pre-image of $G_i$ is of the form 
$\bE_F|_{\cM_3}^{-1}(G_i)=\cM(p_iM_3p_i)$ for $i=1,2$ where $p_i\in M_3$ 
is a projection of rank two (\ref{eq:pre-image}). Since  $L(F)=\bE_F(\cM_3)$ 
and since the images of $p_1$ and $p_2$ intersect in a one-dimensional 
subspace of $\bC^3$, we have $G_1\cap G_2\neq\emptyset$. The intersection 
$G_1\cap G_2$ is a face of $G_1$ and of $G_2$. Large faces being ellipses 
or segments, $G_1\cap G_2$ is an extreme point of $G_1$ and $G_2$.
\par
Let $v_i$ be a unit vector which exposes $G_i$, $i=1,2$. By this we mean,
in the notation of (\ref{eq:expface}), that $G_i=\bF_{L(F)}(v_i)$. Since 
$G_1\cap G_2\neq\emptyset$ and $G_1\neq G_2$, the vectors $v_1$ and $v_2$ 
span a two-dimensional subspace $U\subset\bR^3$. We choose an orthogonal 
transformation $\alpha$, defined in (\ref{eq:trafo2a}), which rotates $U$ 
into the $x_1$-$x_2$-plane and we put $G'_i:=\alpha(G_i)$, $i=1,2$. Using 
the map $\beta$ corresponding to $\alpha$, defined in (\ref{eq:trafo2}), 
we put $F'=(F'_1,F'_2,F'_3):=\beta(F)$. Then $G'_1$ and $G'_2$ are distinct 
intersecting large faces of $L(F')$. Since $G'_1$ and $G'_2$ are exposed 
by vectors in the $x_1$-$x_2$-plane, the projected large faces $P(G'_1)$ 
and $P(G'_2)$ are intersecting proper exposed faces of 
$L(F'_1,F'_2)=P(L(F'))$. Since $G'_i=P|_{L(F')}^{-1}(P(G'_i))$ for $i=1,2$,
we notice that $P(G'_1)\subset P(G'_2)$ implies $G'_1\subset G'_2$. But 
$G'_1\subset G'_2$ is impossible for distinct large faces $G'_1\neq G'_2$,
so $P(G'_1)$ is no singleton and $P(G'_1)\neq P(G'_2)$. Similarly, 
$P(G'_2)\subset P(G'_1)$ is impossible, from which the claim follows.
\hspace*{\fill}$\square$\\
\par
Next we study $L(F)$ without corner points. Mutually distinct large faces 
$G_1,G_2,G_3$ of $L(F)$ satisfy the assumptions of Lemma~\ref{lem:class2}.
\begin{lem}\label{lem:class2}
Let $C\subset\bR^3$ be a convex subset without corner point. Let 
$G_1,G_2,G_3$ be proper exposed faces of $C$, none of which is 
included in any of the others. Then $G_1\cap G_2\cap G_3=\emptyset$.
\end{lem}
{\em Proof.} By contradiction, we assume that $G=G_1\cap G_2\cap G_3$ is 
non-empty. Since $L(F)$ has no corner points, the normal cone of $G$ has
non-maximal dimension $\dim(N(G))\leq 2$. 
As $G$ is strictly included in $G_i$ for $i=1,2,3$, the antitone lattice 
isomorphism (\ref{eq:exp-norm-iso}) shows that $N(G_i)$ is strictly 
included in $N(G)$. Proposition 4.8 of \cite{Weis2012a} shows that $N(G_i)$ 
is a proper face of $N(G)$, so $\dim(N(G_i))<\dim(N(G))\leq 2$ holds for 
$i=1,2,3$. Since the $G_i$ are proper exposed faces of $L(F)$ we have 
$\dim(N(G_i))\geq 1$. Summarizing the dimension count, we have 
$\dim(N(G_i))=1$ for $i=1,2,3$ and $\dim(N(G))=2$. But this is a 
contradiction, as a two-dimensional convex cone cannot have three 
one-dimensional faces.
\hspace*{\fill}$\square$\\
\par
A complete graph with vertex set $\cL(F)$ can be embedded into the 
relative boundary of $L(F)$. 
\begin{lem}[Graph embedding]\label{lem:graph-embedding}
Let $F\in(\her_3)^3$, let $L(F)$ have no corner point, and let $k$ be the 
number of large faces of $L(F)$. Then the complete graph on $k$ vertices 
embeds into the union of large faces with one vertex at the centroid of 
each large face. 
\end{lem}
{\em Proof:} 
For each $G\in\cL(F)$ we denote by $c(G)$ the centroid of $G$ and take 
$c(G)$ as a vertex of the graph to be embedded. Let $G,H\in\cL(F)$ be
distinct. We would like to embed the edge $\{c(G),c(H)\}$ connecting
$c(G)$ and $c(H)$ into the union of large faces as the curve which is 
the union of two segments
\[
\overline{c(G)c(H)}:=[c(G),p(G,H)]\,\cup\,[p(G,H),c(H)]
\]
where $p(G,H)$ is the unique intersection point of $G$ and $H$ found 
in Lemma~\ref{lem:class1}. It remains to show for any 
$G_1,G_2,H_1,H_2\in\cL(F)$ with $G_1\neq H_1$ and $G_2\neq H_2$ and 
with $\{G_1,H_1\}\neq \{G_2,H_2\}$ that the curves 
$\overline{c(G_1)c(H_1)}$ and $\overline{c(G_2)c(H_2)}$ have no 
intersection, except possibly at their end points. Otherwise, by 
construction of the curves, we have $p(G_1,H_1)=p(G_2,H_2)$. Now 
Lemma~\ref{lem:class2} shows $\{G_1,H_1\}=\{G_2,H_2\}$ which completes 
the proof.
\hspace*{\fill}$\square$\\
\par
Each two segments in $\cL(F)$ produce a corner point of $L(F)$ at 
their intersection.
\begin{lem}\label{lem:exclude-two-segments}
Let $F\in(\her_3)^3$ and let there be two distinct segments in $\cL(F)$.
Then the two segments intersect in a corner point of $L(F)$.
\end{lem}
{\em Proof:} 
Lemma~\ref{lem:class1} proves that, after applying an affine transformation,
if necessary, the two segments in $\cL(F)$ project onto the $x_1$-$x_2$-plane 
to two one-dimensional faces of $L(F_1,F_2)$. The classification of the 
numerical range of a $3$-by-$3$ matrix \cite{Kippenhahn1951,Keeler-etal1997} 
shows that either $L(F_1,F_2)$ is a triangle or the convex hull of an ellipse 
and a point outside the ellipse. 
\par
First, let $L(F_1,F_2)$ be a triangle. Another affine transformation allows 
us to take
\[
F_1=\diag(0,0,-1)
\quad\mbox{and}\quad 
F_2=\diag(0,-1,0)
\] 
where the triangle $L(F_1,F_2)$ has vertices $(0,0)$, $(0,-1)$, and $(-1,0)$.
Our strategy is to describe $F_3$ based on the assumption that 
the pre-image $\bF_{L(F)}(1,0,0)$ resp.~$\bF_{L(F)}(0,1,0)$ of the segment
$[(0,0),(0,-1)]$ resp.~$[(-1,0),(0,0)]$ under $P|_{L(F)}$ is a segment 
itself. Let
\[
v:=(1,0,0),
\qquad 
v_1:=(0,1,0), 
\quad\mbox{and}\quad 
v_2:=(0,0,1).
\]
Another way of saying that $\bF_{L(F_1,F_2)}(1,0)=P\,\bF_{L(F)}(1,0,0)$ 
is a segment is to say that the vectors $v,v_1$ span the eigenspace of 
$F_1$ corresponding to the largest eigenvalue and that
\[
\langle v,F_2 v\rangle \neq \langle v_1,F_2 v_1\rangle.
\]
By assumption, 
$\bF_{L(F)}(1,0,0)=P|_{L(F)}^{-1}\,\bF_{L(F_1,F_2)}(1,0)$ is a segment, 
hence (\ref{eq:pre-image}) shows
\[
F_3|_{{\rm span}\{v,v_1\}} \in
{\rm span}\{ F_2|_{{\rm span}\{v,v_1\}}, \id_3|_{{\rm span}\{v,v_1\}} \},
\]
where $A|_X$ denotes the compression to a subspace $X\subset\bC^3$ of 
the linear map defined by $A\in M_3$ in the standard basis. In particular,
\[
\langle v_1,F_3 v\rangle \in
{\rm span}\{ \langle v_1,F_2 v\rangle, \langle v_1,v\rangle \} = \{ 0 \}.
\]
Similarly $\langle v_2,F_3 v\rangle=0$. Since $v_1$ and $v_2$ span the 
orthogonal complement of $v$, we find $F$ in block diagonal form 
\[
F_1=\left(
\begin{array}{c|c}
0 & 0 \hspace{5pt} 0 \\ \hline
0 & \raisebox{-7pt}{{\mbox{{$G_1$}}}} \\[-5pt]
0 & 
\end{array}
\right),
\quad
F_2=\left(
\begin{array}{c|c}
0 & 0 \hspace{5pt} 0 \\ \hline
0 & \raisebox{-7pt}{{\mbox{{$G_2$}}}} \\[-5pt]
0 & 
\end{array}
\right),
\quad
F_3=\left(
\begin{array}{c|c}
a & 0 \hspace{5pt} 0 \\ \hline
0 & \raisebox{-7pt}{{\mbox{{$G_3$}}}} \\[-5pt]
0 & 
\end{array}
\right)
\]
for some $a\in\bR$ and $G=(G_1,G_2,G_3)\in(\her_2)^3$. Since $L(G)$ 
projects to the segment 
\[
P(L(G))=L(G_1,G_2)=[(-1,0),(0,-1)]
\] 
in the $x_1$-$x_2$-plane and because $L(F)$ is the convex hull of 
$L(G)\cup\{(0,0,a)\}$, it follows that $(0,0,a)$ is a corner point 
of $L(F)$.
\par
Second, let $L(F_1,F_2)$ be the convex hull of an ellipse and a point 
outside the ellipse. We have to distinguish a family of affinely 
inequivalent numerical ranges. Lemma 5.1 of \cite{SpitkovskyWeis2015} 
proves that there is a real $b>1$ such that $(F_1,F_2)$ is equivalent 
modulo transformations (\ref{eq:trafo1}) and (\ref{eq:trafo2}) to 
$(F_1',F_2')$, where 
\[
F_1'=\left(\begin{array}{ccc}
b & 0 & 0\\
0 & 0 & 1\\
0 & 1 & 0
\end{array}\right)
\quad\mbox{and}\quad 
F_2'=\left(\begin{array}{ccc}
0 & 0 & 0\\
0 & 0 & -\ii\\
0 & \ii & 0
\end{array}\right).
\]
Another affine transformation allows us to take $F_1$ and $F_2$ 
of the form
\[
F_1=F_1'+\sqrt{b^2-1}\cdot F_2'
\quad\mbox{and}\quad
F_2=F_1'-\sqrt{b^2-1}\cdot F_2'.
\]
Following the same strategy as in the first case of a triangle, we put
\[
v:=(1,0,0),
\qquad
v_1:=(0,1-\ii\sqrt{b^2-1},b),
\quad\mbox{and}\quad
v_2:=(0,1+\ii\sqrt{b^2-1},b).
\]
We obtain $\langle v_1,F_3 v\rangle=0$ because $v$ and $v_1$ span the 
eigenspace of $F_1$ corresponding to the maximal eigenvalue $b>1$ (the 
other eigenvalue of $F_1$ is $-b$) while 
\[
\langle v,F_2 v\rangle 
 = b > 2/b - b 
 = \langle \tfrac{v_1}{\|v_1\|},F_2 \tfrac{v_1}{\|v_1\|}\rangle.
\]
Similarly $\langle v_2,F_3 v\rangle=0$ follows and shows that $F$ 
has block diagonal form 
\[
F_1=\left(
\begin{array}{c|c}
b & 0 \hspace{5pt} 0 \\ \hline
0 & \raisebox{-7pt}{{\mbox{{$G_1$}}}} \\[-5pt]
0 & 
\end{array}
\right),
\quad
F_2=\left(
\begin{array}{c|c}
b & 0 \hspace{5pt} 0 \\ \hline
0 & \raisebox{-7pt}{{\mbox{{$G_2$}}}} \\[-5pt]
0 & 
\end{array}
\right),
\quad
F_3=\left(
\begin{array}{c|c}
a & 0 \hspace{5pt} 0 \\ \hline
0 & \raisebox{-7pt}{{\mbox{{$G_3$}}}} \\[-5pt]
0 & 
\end{array}
\right)
\]
for some $a\in\bR$ and $G=(G_1,G_2,G_3)\in(\her_2)^3$. By construction,
$L(G_1,G_2)$ is a non-degenerate ellipse (tightly) bounded
inside the square $[-b,b]\times[-b,b]$. The numerical range is the 
convex hull of the ellipse and $(b,b)$. Since $(b,b)$ 
lies outside the ellipse, we conclude that $(b,b,a)$ is a corner point 
of $L(F)$.
\hspace*{\fill}$\square$\\
%
%
\section{Examples}
\label{sec:examples}
\par
We analyze three-dimensional examples of joint algebraic numerical ranges
without corner points, one for each configuration of Figure~\ref{fig:classes}. 
Some of the examples are depicted using a heuristic algebraic drawing procedure. 
\par
\begin{figure}
\begin{picture}(13,6)
\put(0,0){\includegraphics[width=6cm]{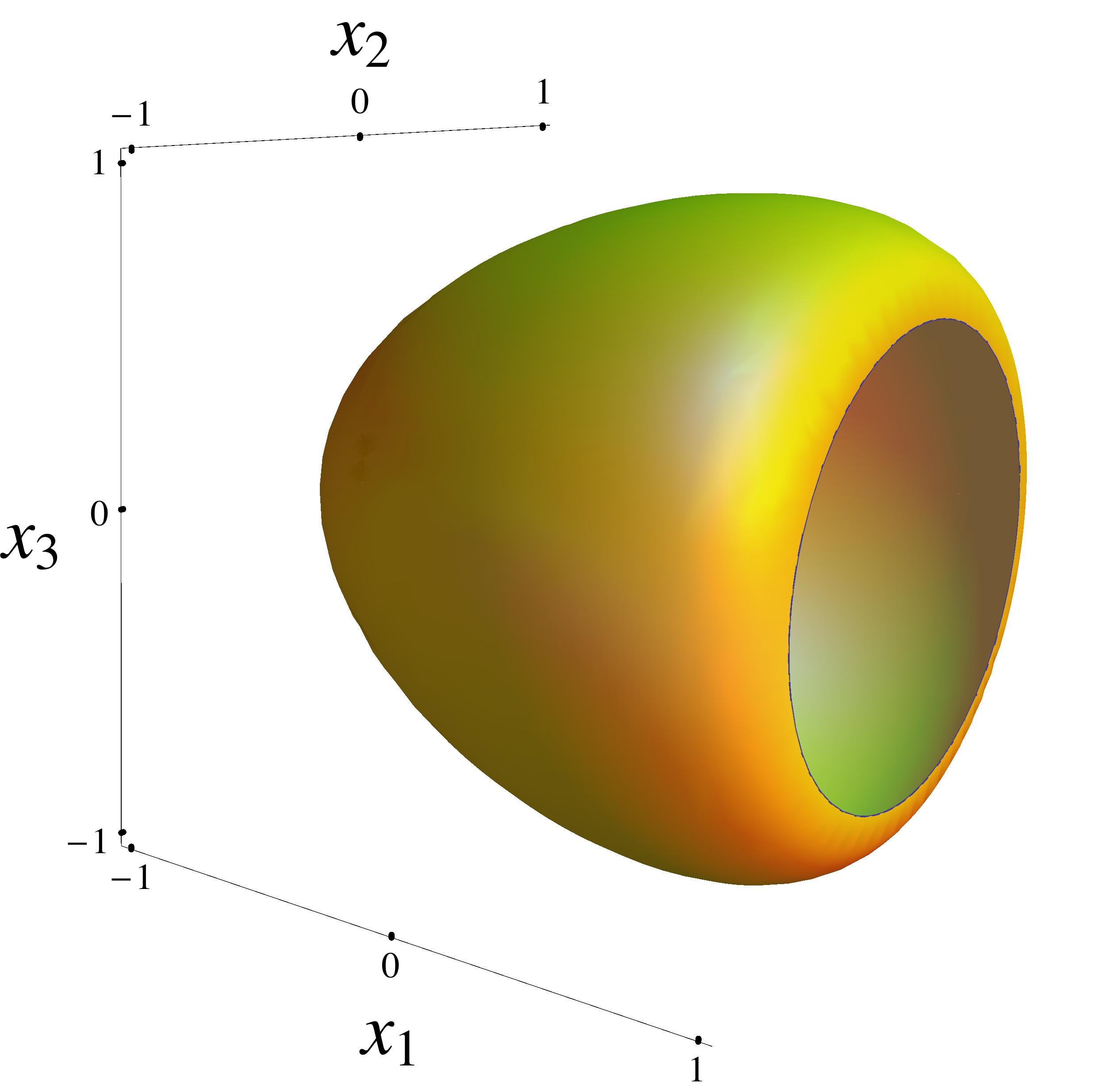}}
\put(0,0){a)}
\put(7,0){\includegraphics[width=6cm]{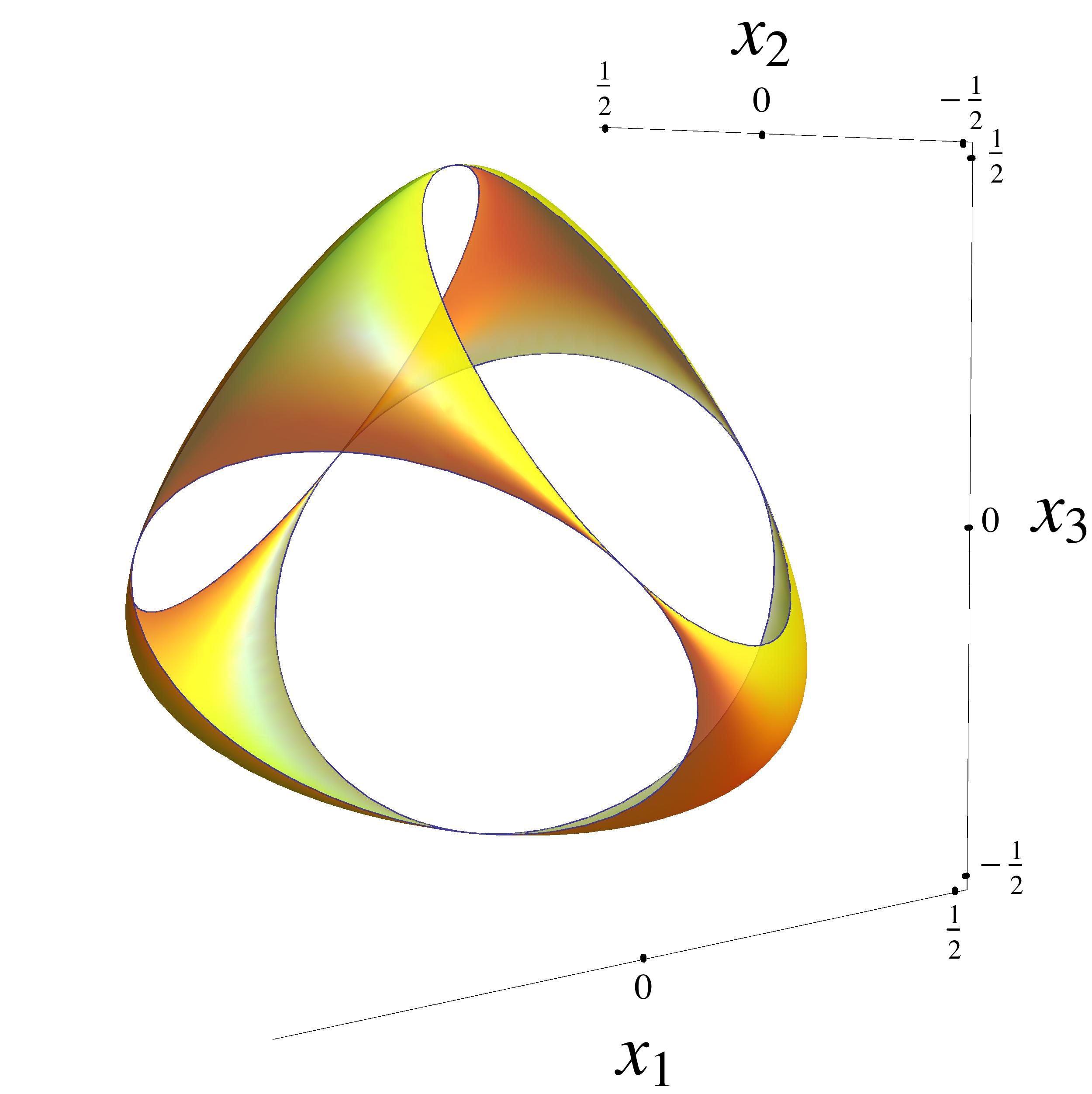}}
\put(7,0){b)}
\end{picture}
\caption{\label{fig:1}%
a) Object with one ellipse and
b) object with four ellipses at the boundary. 
The depicted surfaces are the pieces of the boundary generating
surfaces which lie on the boundary of the joint algebraic numerical 
range.}
\end{figure}
\begin{figure}
\begin{picture}(13,6)
\put(0,0){\includegraphics[height=6cm]{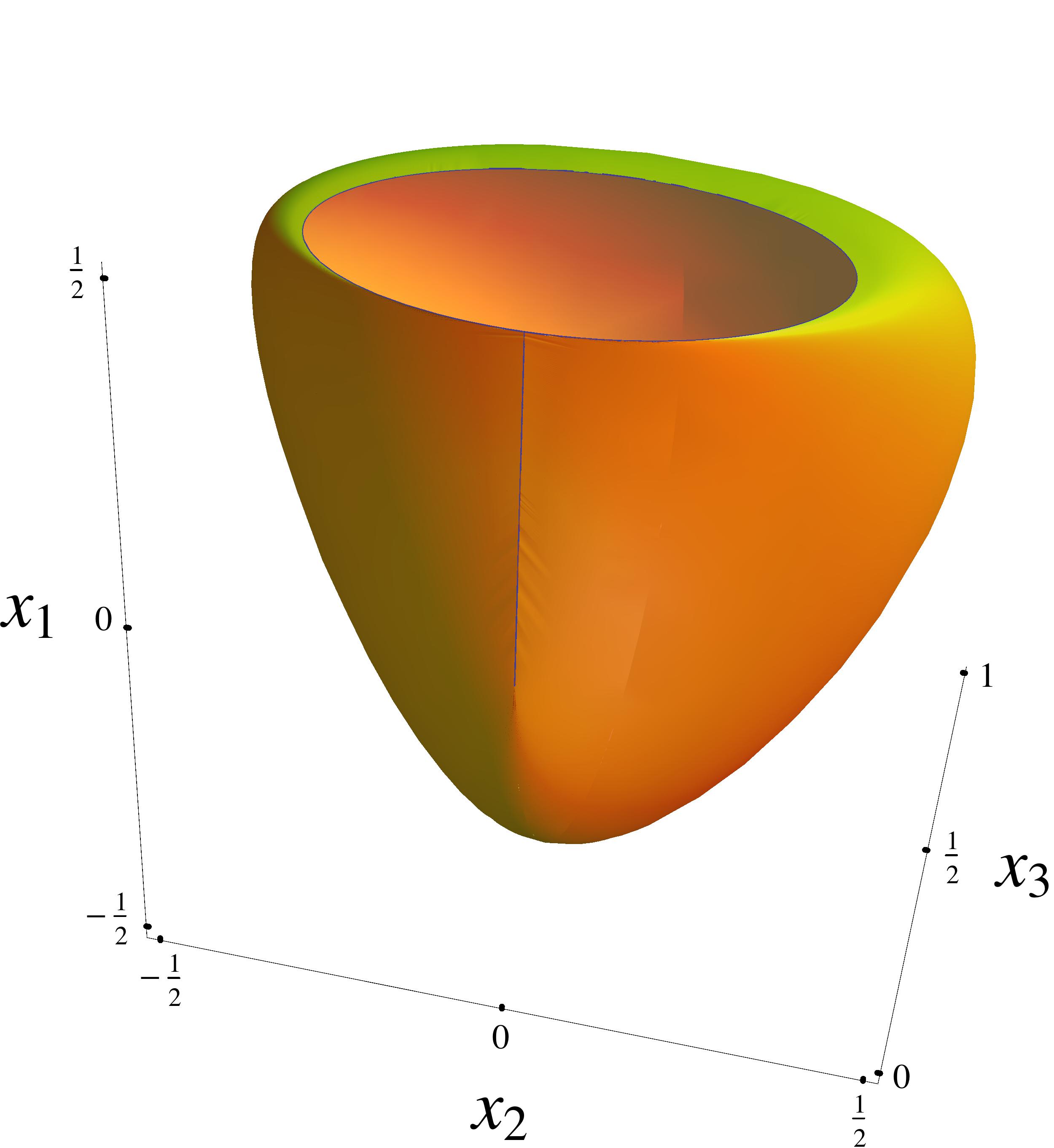}}
\put(0,0){a)}
\put(7,0){\includegraphics[width=6cm]{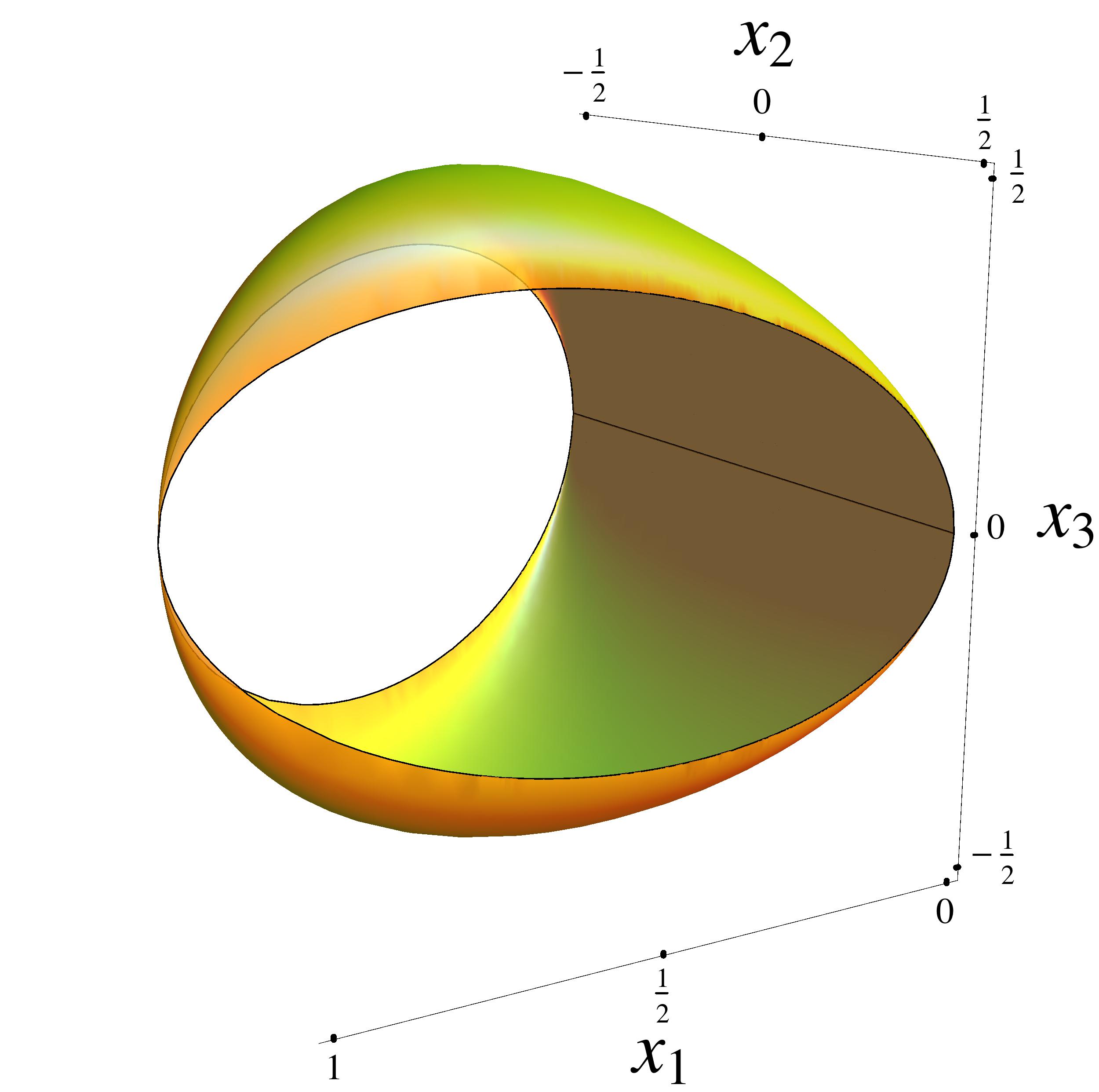}}
\put(7,0){b)}
\end{picture}
\caption{\label{fig:2}%
a) Object with one segment and one ellipse and
b) object with one segment and two ellipses. }
\end{figure}

For all examples we write down the outer normal vectors $u\in\bR^3$ of all 
large faces and we provide hermitian squares as witnesses that there are 
no other large faces. The details are explained in Example~\ref{exa:01} and 
are omitted later on. We also omit the explicit verification that $u$ exposes 
a large face, since this is an easy computation with $2$-by-$2$ matrices 
(\ref{eq:pre-image}):  If the spectral projection of the greatest eigenvalue 
of $F(u)$ is $p$ then the exposed face $\bF_{L(F)}(u)$ is the joint numerical 
range of the compressions of $F_1,F_2,F_3$ to the range of $p$.
\begin{rem}[Heuristic drawing method for joint algebraic numerical ranges]
Recall from Remark~\ref{rem:intro} the definition of the complex projective 
hypersurface $S_F$ in $\bP^3$ with defining polynomial 
$p(u_0,u_1,u_2,u_3):=\det(u_0\id+u_1F_1+u_2F_2+u_3F_3)$, the dual variety 
$S_F^*\subset\bP^n{}^*$, and the {\em boundary generating hypersurface} 
$S_F^*(\bR)\subset\bR^3$. Since $S_F$ is a hypersurface, $S_F^*$ is the
{\em Gauss image} \cite{Harris1995} of $S_F$. We compute a 
Groebner basis\footnote{An algorithm to compute the dual of a variety, which
may not be a hypersurface, is described in \cite{RostalskiSturmfels2012}.} 
of the ideal of polynomials vanishing on $S_F^*$ by eliminating the 
variables $u_0,u_1,u_2,u_3$ from the ideal generated by
\[
p,
\qquad
\partial_{u_i}p-x_i,
\quad
i=0,1,2,3.
\]
In the following examples, the ideal of $S_F^*$ is generated by a
polynomial $\tilde{q}(x_0,x_1,x_2,x_3)$ and the 
{\em boundary generating surface} of $F$ is
\[
S_F^*(\bR)=\{x\in\bR^3\mid q(x)=0\}
\]
where $q(x_1,x_2,x_3):=\tilde{q}(1,x_1,x_2,x_3)$. While for $n=2$ and 
$F\in(\her_d)^2$ the numerical range $L(F_1,F_2)$ is the convex hull
of $S_F^*(\bR)$, the analogue is wrong for $n=3$ because the 
boundary generating surface can contain lines \cite{ChienNakazato2010}. 
The drawings of Figure~\ref{fig:1} and~\ref{fig:2} were generated with
Mathematica. They show pieces of $S_F^*(\bR)$ for which parametrizations 
were obtained. The joint algebraic numerical range $L(F)$ seems to be 
accurately reproduced by the convex hull of these pieces. Pieces of 
$S_F^*(\bR)$ which do not touch the boundary of $L(F)$ were excluded 
from the drawings.
\end{rem}
\par
We discuss examples for the configurations of Figure~\ref{fig:classes} 
with the exception of ovals. Examples of ovals are the Euclidean balls 
in Lemma~\ref{lem:M2-singletons}. More general ovals are discussed in 
\cite{KrupnikSpitkovsky2006}.
%
%
\begin{exa}[$(e,s)=(1,0)$, one ellipse, no segments]
\label{exa:01}
See Figure~\ref{fig:print01} and~\ref{fig:1}a) for 
pictures. If
\[
F_1:=
\left(\begin{smallmatrix}
 1 & 0 & 0 \\
 0 & 1 & 0 \\
 0 & 0 & -1 \\
\end{smallmatrix}\right),
\quad
F_2:=\tfrac{1}{\sqrt{2}}
\left(\begin{smallmatrix}
 0 & 1 & 0 \\
 1 & 0 & 1 \\
 0 & 1 & 0 \\
\end{smallmatrix}\right),
\quad
F_3:=\tfrac{1}{\sqrt{2}}
\left(\begin{smallmatrix}
 0 & -\ii & 0 \\
 \ii & 0 & -\ii \\
 0 & \ii & 0 \\
\end{smallmatrix}\right),
\] 
then
\begin{align*}
q= & -4 x_1^3 - 4 x_1^4 + 27 x_2^2 + 18 x_1 x_2^2 - 13 x_1^2 x_2^2
 - 32 x_2^4 + 27 x_3^2 + 18 x_1 x_3^2 - 13 x_1^2 x_3^2 \\
 & - 64 x_2^2 x_3^2 - 32 x_3^4.
\end{align*}
The greatest eigenvalue of $F_1$ is degenerate and a direct computation 
proves that $\bF_{L(F)}(1,0,0)$ is an ellipse (hence $\bF_{L(F)}(-1,0,0)$ 
is a singleton). The sum of squares representation (\ref{eq:Ilyushechkin}) 
of the modulus $|\delta(F(u))|$ of the discriminant of $F(u)$ contains 
the term 
\[
|M_\nu|^2=(u_2^2 + u_3^2)^3/8
\]
corresponding to 
$\nu:=\{{\scriptstyle (1,1)},{\scriptstyle(1,2)},{\scriptstyle(1,3)}\}$.
This term vanishes only for $u_2=u_3=0$. Thus Remark~\ref{rem:discrim-method} 
shows that $\bF_{L(F)}(1,0,0)$ is the only large face of $L(F)$. 
\end{exa}
%
%
\begin{exa}[$(e,s)=(2,0)$, two ellipses, no segments]
\label{exa:02}
If 
\[
F_1:=
\left(\begin{smallmatrix}
 1 & 0 & 0 \\
 0 & -1 & 0 \\
 0 & 0 & 0 \\
\end{smallmatrix}\right),
\quad
F_2:=
\left(\begin{smallmatrix}
 0 & 0 & 1 \\
 0 & 0 & 0 \\
 1 & 0 & 0 \\
\end{smallmatrix}\right),
\quad
F_3:=
\left(\begin{smallmatrix}
 0 & 1 & 0 \\
 1 & 0 & 0 \\
 0 & 0 & 0 \\
\end{smallmatrix}\right),
\] 
then
\begin{align*}
q= 4 x_1 x_2^2 - 4 x_1^2 x_2^2 - x_2^4 + 4 x_3^2 - 4 x_1^2 x_3^2
 - 4 x_2^2 x_3^2 - 4 x_3^4.
\end{align*}
The hermitian squares corresponding to 
$\nu_1:=\{{\scriptstyle(1, 1)},{\scriptstyle(1, 2)},{\scriptstyle(3, 3)}\}$,
$\nu_2:=\{{\scriptstyle(1, 1)},{\scriptstyle(1, 3)},{\scriptstyle(2, 2)}\}$, 
and
$\nu_3:=\{{\scriptstyle(1, 1)},{\scriptstyle(2, 2)},{\scriptstyle(3, 3)}\}$
are 
\begin{align*}
|M_{\nu_1}|^2 & =(1+u_1^2)^2, && \mbox{if $u_3=\pm1$},\\
|M_{\nu_2}|^2 & =u_2^2 (u_2^2 - 2 u_1^2)^2,\\
\mbox{and} \quad |M_{\nu_3}|^2 & =u_1^2 (u_2^2 - 2 u_1^2)^2, &&\mbox{if $u_3=0$}.
\end{align*}
Thus, $\bF_{L(F)}(-1,\pm\sqrt{2},0)$ are the unique large faces of $L(F)$. 
The $x_1$-axis lies in the boundary generating surface $S_F^*(\bR)$.
\end{exa}
%
%
\begin{exa}[$(e,s)=(3,0)$, three ellipses, no segments]
\label{exa:03}
See Figure~\ref{fig:print03} for a picture corresponding to the matrices 
(\ref{eq:03-b}). If
\[
F_1:=
\left(\begin{smallmatrix}
 0 & 1 & 0 \\
 1 & 0 & 0 \\
 0 & 0 & 0 \\
\end{smallmatrix}\right),
\quad
F_2:=\tfrac{1}{\sqrt{2}}
\left(\begin{smallmatrix}
 0 & \ii & 1\\
 -\ii & 0 & 0\\
 1 & 0 & 0 \\
\end{smallmatrix}\right),
\quad
F_3:=
\left(\begin{smallmatrix}
 1 & 0 & 0 \\
 0 & 1 & 0 \\
 0 & 0 & -1 \\
\end{smallmatrix}\right),
\] 
then
\begin{align*}
q = & -4 x_1^6 - 24 x_1^4 x_2^2 + 27 x_2^4 - 48 x_1^2 x_2^4 
 - 32 x_2^6 + 36 x_1^2 x_2^2 x_3 + 18 x_2^4 x_3 + 8 x_1^4 x_3^2 \\
  & - 4 x_1^2 x_2^2 x_3^2 - 13 x_2^4 x_3^2 - 4 x_2^2 x_3^3 
 - 4 x_1^2 x_3^4 - 4 x_2^2 x_3^4.
\end{align*}
The hermitian squares corresponding to 
$\nu_1:=\{{\scriptstyle(1, 1)},{\scriptstyle(1, 2)},{\scriptstyle(2, 2)}\}$
and
$\nu_2:=\{{\scriptstyle(1, 1)},{\scriptstyle(1, 2)},{\scriptstyle(3, 3)}\}$
are 
\begin{align*}
|M_{\nu_1}|^2 
 & = (1 + 2 u_1^2)/8, 
 && \mbox{if $u_2=\pm1$},\\
\mbox{and} \quad |M_{\nu_2}|^2 
 & = u_1^2(u_1^2 - 4 u_3^2)^2, 
 && \mbox{if $u_2=0$}.
\end{align*}
Thus, $\bF_{L(F)}(0,0,1)$, and $\bF_{L(F)}(\pm2,0,-1)$ are the unique
large faces of $L(F)$. The $x_3$-axis lies in the boundary generating 
surface $S_F^*(\bR)$. Out of curiosity we mention the example
\begin{equation}\label{eq:03-b}\textstyle
F_1:=
\left(\begin{smallmatrix}
 1 & 0 & 0 \\
 0 & 0 & 1 \\
 0 & 1 & 0 \\
\end{smallmatrix}\right),
\quad
F_2:=
\left(\begin{smallmatrix}
 0 & 1 & 0\\
 1 & 0 & 0\\
 0 & 0 & 1 \\
\end{smallmatrix}\right),
\quad
F_3:=
\left(\begin{smallmatrix}
 0 & 0 & \ii \\
 0 & 1 & 0 \\
 -\ii & 0 & 0 \\
\end{smallmatrix}\right).
\end{equation}
Here the normal vectors of the three ellipses are mutually orthogonal 
and $q$ is of degree six with the maximal number of $84$ monomials. 
\end{exa}
%
%
%
\begin{exa}[$(e,s)=(4,0)$, four ellipses, no segments]
\label{exa:04}
See Figure~\ref{fig:print04} and~\ref{fig:1}b) for pictures.
If
\[
F_1:=\tfrac{1}{2}
\left(\begin{smallmatrix}
 0 & 1 & 0 \\
 1 & 0 & 0 \\
 0 & 0 & 0 \\
\end{smallmatrix}\right),
\quad
F_2:=\tfrac{1}{2}
\left(\begin{smallmatrix}
 0 & 0 & 1\\
 0 & 0 & 0\\
 1 & 0 & 0 \\
\end{smallmatrix}\right),
\quad
F_3:=\tfrac{1}{2}
\left(\begin{smallmatrix}
 0 & 0 & 0 \\
 0 & 0 & 1 \\
 0 & 1 & 0 \\
\end{smallmatrix}\right),
\] 
then
\[
q = x_1 x_2 x_3 - x_1^2 x_2^2 - x_1^2 x_3^2 - x_2^2 x_3^2.
\]
For all unit vectors $u\in\bR^3$ the discriminant of $F(u)$ is 
\begin{align*}
\delta(F(u)) = & \tfrac{1}{32} 
((u_1^2 - u_2^2)^2 + (u_2^2 - u_3^2)^2 + (u_3^2 - u_1^2)^2\\
& +6(u_3^2 (u_1^2 - u_2^2)^2 + u_1^2 (u_2^2 - u_3^2)^2 
+ u_2^2 (u_3^2 - u_1^2)^2)). 
\end{align*}
Thus (\ref{eq:pre-image}) proves that 
\begin{align*}
& \bF_{L(F)}(-1,-1,-1), & \bF_{L(F)}(-1,1,1),\\
& \bF_{L(F)}(1,-1,1),  & \mbox{and}\quad \bF_{L(F)}(1,1,-1)
\end{align*}
are the unique large faces of $L(F)$. The boundary generating surface 
$S_F^*(\bR)$ is known as the {\em Roman surface}. It contains the three
coordinate axes. This example was discussed in \cite{Henrion2011}.
\end{exa}
%
%
\begin{exa}[$(e,s)=(0,1)$, no ellipses, one segment]
\label{exa:10}
See Figure~\ref{fig:print10} for a picture. If
\[
F_1:=
\left(\begin{smallmatrix}
 0 & 1 & 0 \\
 1 & 0 & 0 \\
 0 & 0 & 0 \\
\end{smallmatrix}\right),
\quad
F_2:=
\left(\begin{smallmatrix}
 1 & 0 & 0 \\
 0 & -1 & 0 \\
 0 & 0 & 1 \\
\end{smallmatrix}\right),
\quad
F_3:=\tfrac{1}{\sqrt{2}}
\left(\begin{smallmatrix}
 0 & 0 & \ii \\
 0 & 0 & 1 \\
 -\ii & 1 & 0 \\
\end{smallmatrix}\right),
\] 
then $q$ has degree eight and 31 monomials. The hermitian squares 
corresponding to 
$\nu_1:=\{{\scriptstyle(1, 1)},{\scriptstyle(1, 2)},{\scriptstyle(1, 3)}\}$
and
$\nu_2:=\{{\scriptstyle(1, 1)},{\scriptstyle(1, 2)},{\scriptstyle(3, 3)}\}$
are 
\begin{align*}
|M_{\nu_1}|^2 
 & = (1 + 4 u_1^4 + 4 u_1^2 (1 + 4 u_2^2))/8, 
 && \mbox{if $u_3=\pm1$},\\
\mbox{and} \quad |M_{\nu_2}|^2 
 & = u_1^6, 
 && \mbox{if $u_3=0$}.
\end{align*}
Thus, $\bF_{L(F)}(0,1,0)$ is the unique large face of $L(F)$. The
affine hull of $\bF_{L(F)}(0,1,0)$ (which is the line 
$\{x\in\bR^3\mid x_1=0, x_2=1\}$) and the $x_2$-axis lie in the 
boundary generating surface $S_F^*(\bR)$. 
\end{exa}
%
%
%
%
\begin{exa}
[$(e,s)=(1,1)$, one ellipse, one segment]
\label{exa:11}
If $\lambda\in\bR$ and 
\[
F_1:=\tfrac{1}{2}
\left(\begin{smallmatrix}
 \lambda & 0 & 0 \\
 0 & 0 & 1 \\
 0 & 1 & 0 \\
\end{smallmatrix}\right),
\quad
F_2:=\tfrac{1}{2}
\left(\begin{smallmatrix}
 0 & 0 & 1 \\
 0 & 0 & 0 \\
 1 & 0 & 0 \\
\end{smallmatrix}\right),
\quad
F_3:=
\left(\begin{smallmatrix}
 0 & 0 & 0 \\
 0 & 0 & 0 \\
 0 & 0 & 1 \\
\end{smallmatrix}\right),
\] 
then
\[
q=-4 x_1^2 x_3^2 - 4 x_2^2 x_3^2 + 4 x_3^3 - 4 x_3^4 + 
 4 x_1 x_2^2 x_3 \lambda - x_2^4 \lambda^2.
\] 
See Figure~\ref{fig:print11} and~\ref{fig:2}a) for pictures, the
latter at $\lambda=1$. For $\lambda=0$, equation~(\ref{eq:vectorsD3real})
of Lemma~\ref{lem:M2-singletons} shows that $L(F)$ is the Euclidean ball 
of radius $\tfrac{1}{2}$ centered at $(0,0,\tfrac{1}{2})$. For $\lambda=1$ 
the hermitian squares corresponding to 
$\nu_1:=\{{\scriptstyle(1, 1)},{\scriptstyle(1, 2)},{\scriptstyle(1, 3)}\}$,
$\nu_2:=\{{\scriptstyle(1, 1)},{\scriptstyle(1, 3)},{\scriptstyle(2, 2)}\}$, 
and
$\nu_3:=\{{\scriptstyle(1, 1)},{\scriptstyle(2, 2)},{\scriptstyle(2, 3)}\}$
are 
\begin{align*}
|M_{\nu_1}|^2 & = u_1^2/64, && \mbox{if $u_2=\pm1$},\\
|M_{\nu_2}|^2 & = 1/64, && \mbox{if $u_2=\pm1$, $u_1=0$},\\
\mbox{and} \quad |M_{\nu_3}|^2 & =u_1^4 u_3^2/16, &&\mbox{if $u_2=0$}.
\end{align*}
Thus, $\bF_{L(F)}(1,0,0)$ and $\bF_{L(F)}(0,0,-1)$ are the unique large 
faces of $L(F)$. The $x_1$-axis lies in the boundary generating surface 
$S_F^*(\bR)$. For $\lambda=2$, the joint algebraic numerical range $L(F)$ 
is affinely isomorphic to the example in Section~3 of 
\cite{ChienNakazato2010}, where the first line in $S_F^*(\bR)$ was 
discovered.
\end{exa}
%
%
\begin{exa}
[$(e,s)=(2,1)$, two ellipses, one segment]
\label{exa:12}
See Figure~\ref{fig:print12} and~\ref{fig:2}b)
for pictures. If 
\[
F_1:=\left(\begin{smallmatrix}
0 & 0 & 0\\
0 & 0 & 0\\
0 & 0 & 1
\end{smallmatrix}\right),
\quad
F_2:=\tfrac{1}{2}\left(\begin{smallmatrix}
0 & 1 & 0\\
1 & 0 & 0\\
0 & 0 & 0
\end{smallmatrix}\right),
\quad
F_3:=\tfrac{1}{2}\left(\begin{smallmatrix}
0 & 0 & 1\\
0 & 0 & 0\\
1 & 0 & 0
\end{smallmatrix}\right),
\]
then
\[
q = -x_1^2 x_2^2 + x_1 x_3^2 - x_1^2 x_3^2 - x_3^4.
\]
The hermitian squares corresponding to 
$\nu_1:=\{{\scriptstyle(1, 1)},{\scriptstyle(1, 3)},{\scriptstyle(2, 2)}\}$
and
$\nu_2:=\{{\scriptstyle(1, 1)},{\scriptstyle(1, 2)},{\scriptstyle(3, 3)}\}$
are 
\begin{align*}
|M_{\nu_1}|^2 
 & = 1/64, 
 && \mbox{if $u_3=\pm1$},\\
\mbox{and} \quad |M_{\nu_2}|^2 
 & = u_2^2(u_2^2-4 u_1^2)^2/64, 
 && \mbox{if $u_3=0$}.
\end{align*}
Thus, $\bF_{L(F)}(-1,0,0)$ and $\bF_{L(F)}(1,\pm2,0)$ are the unique
large faces of $L(F)$. The $x_1$- and $x_2$-axes lie in the boundary 
generating surface $S_F^*(\bR)$. The joint algebraic numerical range 
$L(F)$ is affinely isomorphic to the object in Example~6 of 
\cite{Chen-etal2015}.
\end{exa}
%
%
%
\bibliographystyle{plain}

%
%
%
%
\vspace{1cm}
\parbox{8cm}{%
Konrad Szyma\'nski\\
Marian Smoluchowski Institute of Physics\\
Jagiellonian University\\
\L ojasiewicza 11\\
30-348 Kraków\\ 
Poland\\
e-mail: konrad.szymanski@uj.edu.pl}
\vspace{\baselineskip}
\par\noindent
\parbox{10cm}{%
Stephan Weis\\
Centre for Quantum Information and Communication\\
Université libre de Bruxelles\\
50 av.~F.D.~Roosevelt - CP165/59\\
1050 Bruxelles\\
Belgium\\
e-mail: maths@weis-stephan.de}
\vspace{\baselineskip}
\par\noindent
\parbox{8cm}{%
Karol {\.Z}yczkowski\\
Marian Smoluchowski Institute of Physics\\
Jagiellonian University\\
\L ojasiewicza 11\\
30-348 Kraków\\ 
Poland\\
e-mail: karol.zyczkowski@uj.edu.pl\\
~\\
and\\
~\\
Center for Theoretical Physics\\
of the Polish Academy of Sciences\\
Al.\ Lotnikow 32/46\\ 
02-668 Warsaw\\		 
Poland}
\end{document}